\documentclass[11pt,a4paper]{amsart}

\title{Geometric Combinatorics of Polynomials with Rotationally Symmetric Roots}
\author{Benedict Lee}
\usepackage{geometry}
\usepackage{physics}
\usepackage{amsthm}%for theorems
\usepackage{mathtools}
\usepackage{tikz-cd}
\usepackage{amsmath,amsfonts,amssymb}
\usepackage[
backend=biber,
style=alphabetic,
]{biblatex}
\newcommand{\R}{\mathbb{R}}
\newcommand{\C}{\mathbb{C}}
\newcommand{\D}{\partial}

\newcommand{\NCP}{\textsc{NCPart}}
\newcommand{\cvl}{\text{\textbf{cvl}}}
\newcommand{\cpt}{\text{\textbf{cpt}}}
\newcommand{\re}{\textsc{Re}}
\newcommand{\imag}{\textsc{Im}}
\newcommand{\LL}{\textsc{LL}}

\newcommand{\cox}{\textsc{Cox}}
\newcommand{\art}{\textsc{Art}}
\newcommand{\poly}{\textsc{Poly}}
\newcommand{\rshp}{\textsc{RedShape}}
\newcommand{\polyk}[1]{\textsc{Poly}^{mc}_{d,{#1}}}

\newcommand{\multi}{\textsc{Multi}}
\newcommand{\sym}{\textsc{Sym}}

\newcommand{\closedrect}{\raisebox{-0.2em}{\includegraphics[scale=0.2]{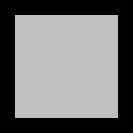}}}
\newcommand{\openrect}{\raisebox{-0.2em}{\includegraphics[scale=0.2]{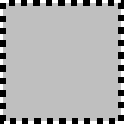}}}
\newcommand{\closedann}{\raisebox{-0.2em}{\includegraphics[scale=0.2]{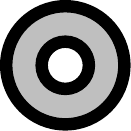}}}

\newcommand{\closedint}{\raisebox{0.1em}{\includegraphics[scale=0.2]{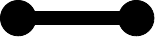}}}
\newcommand{\circcle}{\raisebox{-0.2em}{\includegraphics[scale=0.2]{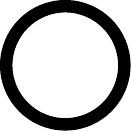}}}

\newcounter{bencomments}

\theoremstyle{definition}
\newtheorem{thm}{Theorem}[section]
\newtheorem{cor}[thm]{Corollary}
\newtheorem{definition}[thm]{Definition}
\newtheorem{prop}[thm]{Proposition}
\newtheorem{lemma}[thm]{Lemma}
\newtheorem{example}[thm]{Example}
\newtheorem{rmk}[thm]{Remark}

\begin{document}

%%%%%%%%%%%%%%%%%
\begin{abstract}
    In a recent paper by Dougherty and McCammond, the authors introduced a piecewise Euclidean cell structure for the space of monic centered polynomials. This was achieved by studying the relationship between the polynomials and their critical values, using the Lyashko-Looijenga ($\LL$) map, as well as the noncrossing combinatorics of the polynomials. In this paper, we extend their results to the space of polynomials with rotationally symmetric roots by studying an appropriate generalisation of the $\LL$ map, which relates the polynomials and their reduced critical values. The combinatorics of the polynomials are noncrossing partitions of type $B_d$. This give us a piecewise Euclidean cell structure built using the order complex of the noncrossing partition of type $B_d$. Along the way, we also give a compactification of this space of polynomials. 
\end{abstract}

\maketitle

%%%%%%%%%%%%%%%%%%%%%%
\section{Introduction}
%%%%%%%%%%%%%%%%%%%%%%
Let $\poly_d^{mc}(\C)$ denote the space of monic centered polynomials of degree $d$ with complex coefficients. Recently, Dougherty and McCammond \cite{DM24} gave a piecewise Euclidean cell complex to this space built out of products of \textit{orthoschemes}, which they called the \textit{branched rectangle complex}. This complex is obtained via two steps: the first step involves studying the relationship of $\poly_d^{mc}(\C)$ with their critical values using a map known as the \textit{Lyashko-Looijenga map}. The second step looks at the combinatorics of each polynomial which gives barycentric coordinates into a product of the order complexes of the noncrossing partition of $d$-elements, $\abs{\NCP_d}$.

The construction in \cite{DM24} allowed the authors to obtained a compactification of the space of polynomials, while the cell structure allowed them to describe a gluing of the cells to give a compactification of the quotient hyperplane complement of the symmetric group. Importantly, the cell structure also allowed the authors to construct a homotopy between the quotient compactified hyperplane complement and the dual braid complex, both of which are well-known classifying spaces of the braid group. 

In \cite{DM24}, the authors conjectured that their construction can be generalised to other finite complex reflection groups. In this paper, we extend the results of \cite{DM24} to finite complex reflection groups of type $G(k,1,d)$. Let us denote $\poly_{d,k}^{mc}(\C)$ to be the space of polynomials of degree $dk$ of the form $p(z^k)$, where $p\in\poly_d^{mc}(\C)$.  Following the authors of \cite{DM24}, we denote $\poly_{dk}^{mc}(U)$ to be the subspace of $\poly_{d,k}^{mc}(\C)$ consisting of polynomials with critical values lying inside the subspace $U\subset\C$. We also adopt their symbols $\closedint$ to be a closed interval $[a,b]\subset \R\subset \C$, $\closedrect$ to be a closed rectangle in $\C$, and $\closedann$ to be a closed annulus in $\C$. We show the following three theorems:

\begin{thm}[Interval Theorem]
	The space $\poly_{d,k}^{mc}(\closedint)$ is homeomorphic to the order complex of the noncrossing partition of type $B_d$, denoted as $\abs{\NCP(B_d)}$.
\end{thm}

\begin{thm}[Rectangular Theorem]
The space $\polyk{2}(\closedrect)$ embeds into the product of the order complexes of the noncrossing partition of type $B_d$,  $\abs{\NCP(B_d)}\times\abs{\NCP(B_d)}$.
\end{thm}
\begin{thm}[Annular Theorem]
The space $\polyk{2}(\closedrect)$ embeds into the product of the dual Artin complex and the order complex, $\mathcal{K}_{B_d}\times \abs{\NCP(B_d)}$, where $\mathcal{K}_{B_d}$ is the dual Artin complex of type $B_d$.
\end{thm}

While the three theorems are closely related to the three theorems found in \cite{DM24}, the results of this paper cannot be deduced directly from their paper, and new constructions have to be made to prove these results. The paper is organised as follows: we first review the work of \cite{DM24} as well as the notation used in their paper. Next, we look at the noncrossing combinatorics of $G(k,1,d)$. We then introduce the map that relates polynomials and their reduced critical values, which allows us to compactify our space of polynomials. The noncrossing combinatorics of polynomials in $\poly_{d,k}^{mc}(\C)$ is studied, which allows us to prove the three theorems stated above. 

\textbf{Acknowledgements.} The author would like to thank his advisor, Jon McCammond, for the numerous insightful conversations they had, as well as for his guidance in writing the author's dissertation and this article.

\textbf{Statement on AI use.} No AI was used in the preparation and writing of this article.

%%%%%%%%%%%%%%%%%%%%%%%%%%%%%%%%%%%%%%%%%%%%%%%%%%%%%%%%%%%%%%%%%
\section{Geometric Combinatorics of Space of Complex Polynomials}
\subsection{Topology of the space of complex polynomials}
%%%%%%%%%%%%%%%%%%%%%%%%%%%%%%%%%%%%%%%%%%%%%%%%%%%%%%%%%%%%%%%%%

In this part, we review the results of \cite{DM24}, as well as their notations. The reader is encouraged to look at their paper for the details of the proofs and the constructions.
\begin{definition}[Space of polynomial]
	Let $U\subset \C$, the space $\poly_d^{mc}(U)$ is the set of all monic polynomial of degree $d$ that are \textit{centered}, i.e. their roots sum to zero, whose critical values lie inside $U$.
\end{definition}
Let $X$ be a space, e.g. $\C$ or subsets of $\C$, a \textit{multiset} of size $n$ is a point in the quotient space $X^n/\sym_n$, where $\sym_n$ is the symmetric group of $n$-elements. The space of all multisets of size $n$ is denoted as $\multi_n(X)$. A multiset, $M$, is written as $\{x_1^{k_1},\cdots,x_m^{k_m}\}$, where $k_i$, known as a \textit{multiplicity}, is the number of times $x_j$ appear in a representative of $M$, and $k_1+\cdots k_m=n$. 

\begin{definition}[Shape]
	Given a multiset of size $n$, its multiplicities form a partition of $n$. An integer partition of $n$, $\lambda\in \textsc{IntPart}_n$, is the \textit{shape} of the multiset if it is equal to the multiplicities of the multiset. Similarly, given a set partition $[\lambda]\vdash [n]$, the size of each block partitions $n$, which is also called the \textit{shape} of the set partition, denoted $\textsc{Shape}([\lambda])=\lambda$.
\end{definition}
\begin{definition}[Order on Integer Partitions] Let $\lambda$ be an integer partition and $\pi$ be a partition of the parts of $\lambda$, then there is a new partition $\mu$ whose parts are the sums of the numbers in the block of $\pi$, we denote this as $\lambda \xrightarrow{\pi}\mu$. We say that $\lambda\leq \mu$ if there exists a partition $\pi$ such that $\lambda\xrightarrow{\pi}\mu$. Note that we sometimes write $\lambda\rightarrow\mu$ if the partition $\pi$ is implicitly understood.
\end{definition}
\begin{definition}{[Order on Set Partition]}
	Given two set partitions $[\lambda],[\mu]\vdash[n]$, we can order them by $[\lambda]\leq[\mu]$ if $[\lambda]$ is a refinement of $[\mu]$, i.e. if $B$ is a block in $[\lambda]$, then $B$ is contained in some block in $[\mu]$.
\end{definition}

The space $\poly_d^{mc}(\C)$ is stratified by the shape of the critical values and the critical points of the polynomials.

\begin{definition}[Stratification]
	The space of polynomials has a double stratification by critical point shape and critical value shape. Let $\cpt(p)$ and $\cvl(p)$ denote the multiset of critical points and critical values respectively, and let $\lambda=\lambda(p)$ and $\mu=\mu(p)$ be the shape of $\cpt(p)$ and $\cvl(p)$ respectively. We denote $\lambda(p)\xrightarrow{p}\mu(p)$ to represent the polynomial $p$ mapping $\cpt(p)$ with shape $\lambda$ to $\cvl(p)$ with shape $\mu$, this gives us an ordering $\lambda \rightarrow\mu$. Let $\poly_{\lambda\to\mu}(\C)$ denote the set $\{p\in\poly_{d}\ \vert\ \lambda\to\mu \text{ is }\lambda(p)\xrightarrow{p}\mu(p)\}$, we get the following stratification of our polynomial space
	\begin{align*}
		\poly_d^{mc}(\C)=\bigsqcup_{\substack{\lambda,\mu\vdash d\\\lambda\rightarrow\mu}}\poly_{\lambda\to\mu}(\C)
	\end{align*}
	
	The space of multiset is also stratified. Let us denote $\multi_\mu(\C)$ to be the subspace of $\multi_n(\C)$ consisting of multisets with shape $\mu$, in other words $\multi_\mu(\C)$ denotes the set $\{M\in\multi_n(\C)\ \vert\ \textsc{Shape}(M)=\mu\}$, then 
	\begin{align*}
		\multi_n(\C)=\bigsqcup_{\mu\vdash n}\multi_\mu(\C)
	\end{align*}
\end{definition}

The critical points and the critical values are related by the Lyashko-Looijenga map:

\begin{definition}
	The Lyashko-Looijenga ($\LL$) map is the map sending each polynomial in $\poly_d^{mc}(\C)$ to its multiset of critical values. In symbols:
	\begin{align*}
		\LL:\poly_d^{mc}(\C)& \to\multi_n(\C)\\
		p&\mapsto \cvl(p)
	\end{align*}
	where $\cvl(p)=\{x_1^{k_1},\cdots,x_m^{k_m}\}$ is a multiset where $x_i$ are the critical values of $p$ and $k_i$ is the multiplicity of the critical value. 
\end{definition}
\begin{rmk}
	The $\LL$ map can be confusing when first encountered. Since a multiset of size $n$ in $\C$ is equivalent to $n$ (counted with multiplcities) unordered points in $\C$, which can be further identified with the roots of a monic degree $n$ polynomial, the $\LL$ map can be thought of as a map from the $(d-1)$-tuple of coefficients of $p$ to the $n$-tuple of coefficients of a monic polynomial whose roots are the critical values of $p$.
\end{rmk}
\begin{thm}[Stratified Covering, Theorem 6.19 \cite{DM24}]\label{thm:LLStratCov}
	Let $\lambda,\mu\vdash d$ such that $\lambda\rightarrow \mu$. The restricted map $\LL:\poly_{\lambda\to\mu}^{mc}(\C)\to\multi_\mu(\C)$ is a covering map.
\end{thm}
\begin{rmk}
	While the $\LL$ map is not a covering map, the fact that it is a stratified covering map allowed \cite{DM24} to construct unique path lifting and homotopy lifting under specific circumstances\footnote{They can be constructed when the paths or homotopies are nonsplitting, see \cite{DM24} for the precise definitions.}. This allowed them to homotope the space $\poly_d^{mc}(\C)$ to $\poly_d^{mc}(\openrect)$. 
\end{rmk}
\begin{thm}[Compactification, Theorem 10.9 \cite{DM24}]\label{thm:compactification}
	$\poly_d^{mc}(\closedrect)$ is a compactification of $\poly_d^{mc}(\C)$.
\end{thm}
\subsection{Geometric combinatorics of complex polynomials}
In this section, we fix a closed rectangle $\closedrect$ inside $\C$. In \cite{DM24}, the authors gave every polynomial in $\poly_d^{mc}(\closedrect)$ two chains of noncrossing partitions. As the construction is quite technical, let us use an example to illustrate how the combinatorics of a polynomial arise, and refer the readers to \cite{DM24} for definitions, notations, and details of the construction.

\begin{example}\label{ex:AnExample}
	Consider the polynomial:
	\begin{align}\label{eqn:RunningEgAn}
		p(z)=z^4+ ( -3.182  -0.432i)z^2+
		(0.899 -0.259i)z+   1.630+0.633i
	\end{align}
	which has following multiset of critical values:
	\begin{align*}
		\text{\textbf{cvl}}(p)\coloneqq\{-2.038+0.217i,  0.277-0.281i,
		1.683+0.589i\}
	\end{align*}
	Let $R$ be a closed rectangle containing $\closedrect$ in its interior. In this example, we chose $R=[-3,3]\times i[-3,3]$. In Figure~\ref{fig:PullbackPartition}, $R$ is subdivided into two rectangles, an orange rectangle, which we call $R_2^L$, and a red rectangle. The preimage of $R$ is the $16$-sided `polygon' $P$, and the preimage of $R_2^L$ partitions the $l_i$'s inside $P$ into path connected blocks. The `sides' of $P$ are labelled based on the side of $R$ that it is sent to.  Since all the branching occurs at the critical values, we define a sequence of $R_i^L$ such that each subsequent rectangle contains new critical value(s). This is repeated for a horizontal subdivision of the rectangle. 
	\begin{figure}
		\centering
		\raisebox{1em}{\includegraphics[width=0.55\textwidth]{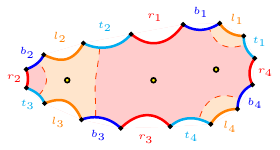}}
		\raisebox{6.5em}{$\xrightarrow{p(z)}$}	\includegraphics[width=0.35\textwidth]{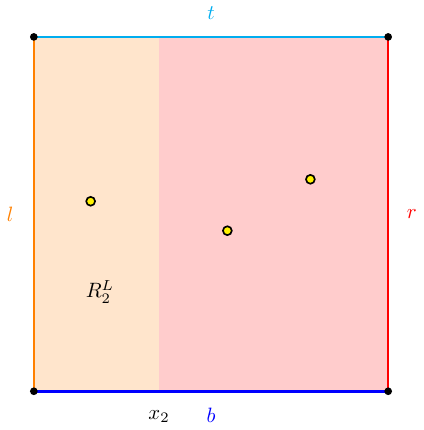}
		\caption{The closed rectangle, $R$, on the right and its preimage polygon, $P=p^{-1}(R)$, on the left. The left subrectangle $R_2^L$ is defined by our choice of subdivision of the interval $[-3,3]$ is shown on the right.}
		\label{fig:PullbackPartition}
	\end{figure}
	
	$R$ is subdivided into left and right subrectangles by partitioning the real and imaginary parts of $R$, i.e. the interval $[-3,3]$. To partition the real part of $R$ to get the left subrectangles, we choose a partition $x_0<\cdots<x_{5}$, such that each subinterval $[x_i,x_{i+1}]$, with $x_i\neq  -3$ or $x_{i+1}\neq 3$, contains only one value of $\re(\cvl(p))$. The \textit{left subrectangle}, $R_i^L,$ is defined to be the rectangle $[-3,x_i]\times i[-3,3]$, with $1\leq i\leq 4$. This process is repeated for the values in $\imag(\cvl(p))$. Figure~\ref{fig:critbox2} shows one such subdivision of $R$. In our example, we choose the following the partitions $x_i$ for the real part and $y_j$ for the imaginary part.
	\begin{figure}
		\centering
		\includegraphics[width=0.5\linewidth]{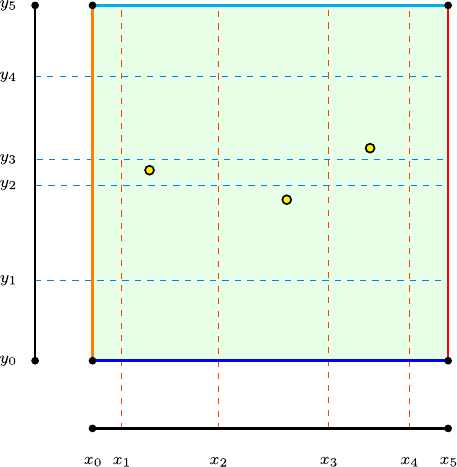}
		\caption{Given critical values lying in $\closedrect$, and a closed rectangle $R=[x_0,x_5]\times i[y_0,y_5]$ containing $\closedrect$ in its interior, we can subdivide $R$ into subrectangles containing different numbers of critical values. More precisely, consider the real interval, $[x_0,x_5]$, and the projection $\re(\cvl(p))$. We subdivide $[x_0,x_5]$ so that the subintervals not containing $x_0$ or $x_5$ contains exactly one value in $\re(\cvl(p))$. The vertical line $\re(z)=x_i$ defines the vertical edge of our subrectangles $R_i^L$.}
		\label{fig:critbox2}
	\end{figure}
	\begin{alignat*}{6}
		x_0&=-3&&<-2.519&&<-0.881&&<0.980&&<2.341&&<  3=x_5\\
		y_0&=-3&&<-1.641&&<-0.032&&<0.403&&<1.794&&<  3=y_5          
	\end{alignat*}
	
	Referring back to Figure~\ref{fig:PullbackPartition}, we see that $l_2$ and $l_3$, the preimages of the left side of $R$, are in the same path connected component of $p^{-1}(R_2^L)$, while $l_1$ and $l_4$ are in distinct components. This is a partition of the $l_i$, which we write as $[\lambda_2]^L=1\vert 23\vert 4$.
	
	Since we are only interested in the combinatorial information, we shall represent the  $P$ in a stylised manner. Figure~\ref{fig:combinedpartitionAn} shows the pull back of $R_i^L$ for the subdivision chosen above, as well as the noncrossing partition obtained from looking at the path connected components.
	
	\begin{figure}
		\centering
		\includegraphics[width=\linewidth]{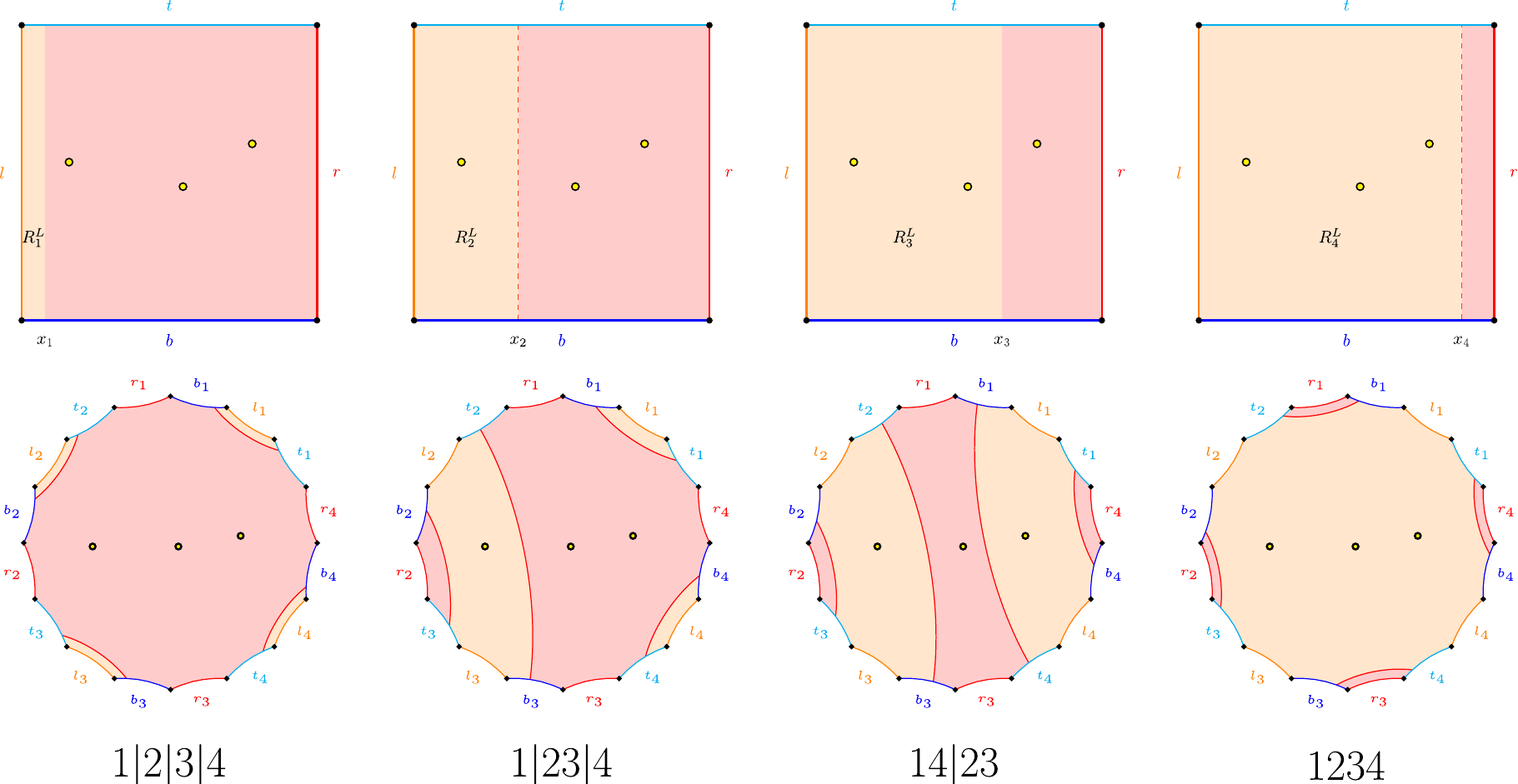}
		\caption{The top row of the diagram shows the various left subrectangles $R_i^L$ from the subdivision chosen earlier. The subrectangles pullback to various polygons which partition $l_i$'s, as shown in the second row. The last row is the noncrossing partition of the $l_i$'s.}
		\label{fig:combinedpartitionAn}
	\end{figure}
	
	The bottom subrectangles $R_j^B$ give a sequence of noncrossing partitions of $b_i$'s as shown in Figure~\ref{fig:combinedpartitionTBAn}.
	
	\begin{figure}
		\centering
		\includegraphics[width=\linewidth]{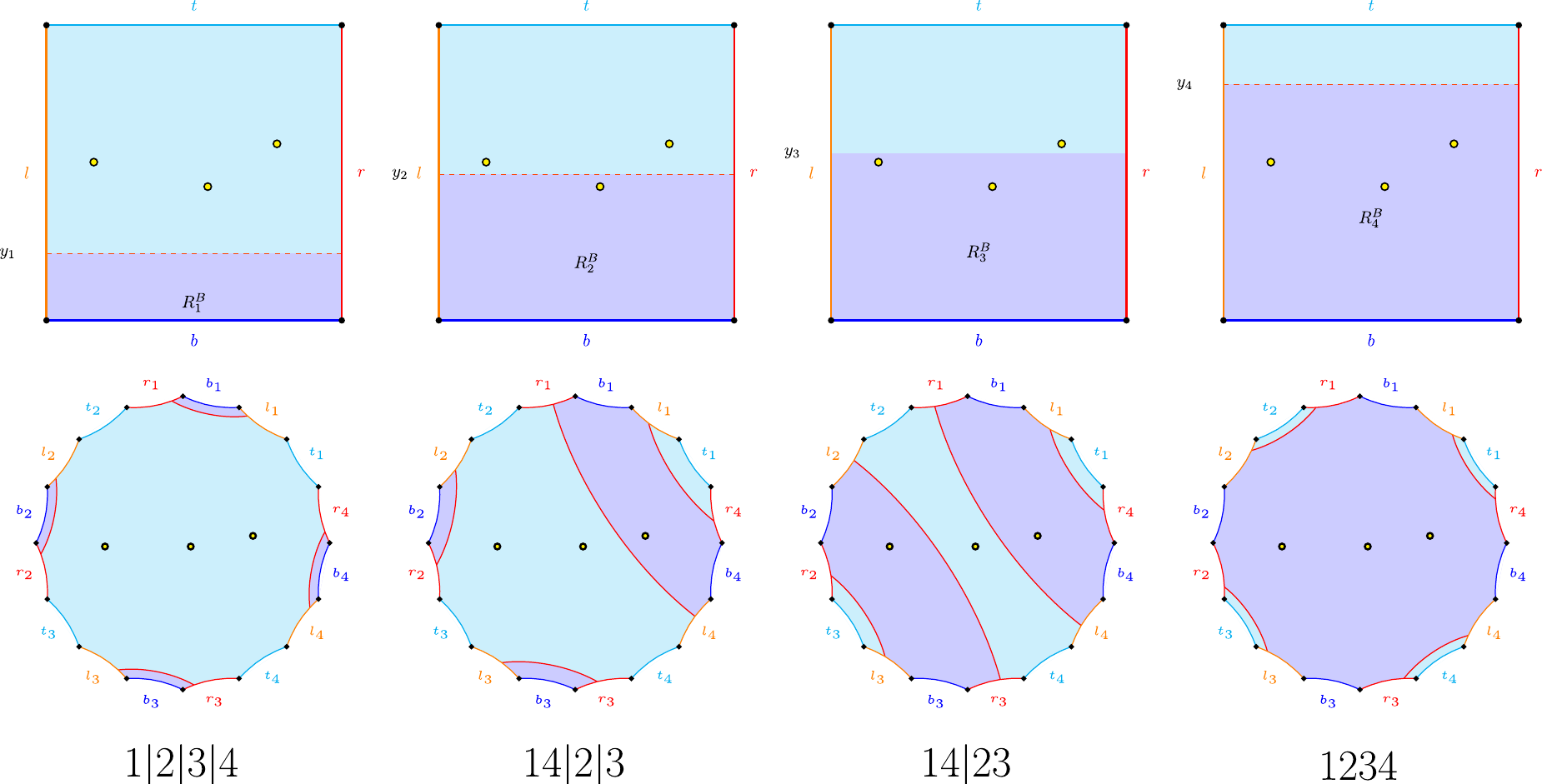}
		\caption{The top row of the diagram shows the various left subrectangles $R_j^B$ from the subdivision chosen earlier. The subrectangles pullback to various polygons which partition $b_j$'s, as shown in the second row. The last row is the noncrossing partition of the $b_j$'s.}
		\label{fig:combinedpartitionTBAn}
	\end{figure}
\end{example}

The chains of noncrossing partitions behave in a regular way:

\begin{prop}\label{prop:strictInc}
	Let $R_1^L\subset R_2^L\subset\cdots\subset R_k^L\subset R$ be a sequence of left subrectangles defined above and let $[\lambda_1]^L,\cdots,[\lambda_k]^L$ be the side partitions obtained, the sequence of partition is noncrossing and strictly increasing, starting from the discrete partition and ending at the indiscrete partition.
\end{prop}

\begin{rmk}[Permutations]\label{rmk:sidePerm}
	Since a complex polynomial $p$ can be turned into a covering map if we remove $\cvl(p)$ from the codomain, and $p^{-1}(\cvl(p))$ from the domain, after a choice of base point $z_a$ and an appropriate numbering of $p^{-1}(z_a)$, the noncrossing partitions can be turned into a sequence of deck transformation given by permutation of $d$-elements. For example, the basepoint $z_a$ may be the bottom left corner of $R$. The boundary curves of our subrectangles are loops based at $z_a$, and they permute the set $p^{-1}(z_a)$. If we label the points in $p^{-1}(z_a)$ such that a preimage is labelled $i$ if it is the boundary of the edges $l_i$ and $b_i$, then the permutations are related to the noncrossing partitions obtained above by replacing the $k$-blocks with $k$-cycle ordered counter-clockwise. Moreover, the sequence of permutations give a factorisation of the $d$-cycle. See Definition 4.17 of \cite{DM24}.
\end{rmk}
\begin{rmk} [Transitive Hurwitz Action and $\LL$ map]
	The existence of stratified path lifting using the $\LL$ map allows one to define a braid group action on the critical values, which gives a braid group action on the factorisation of the $d$-cycle known as the \textit{Hurwitz action}. It is a well-known fact that the Hurwitz action on the \textit{reflection} factorisation of the $d$-cycle is transitive. In particular, there are $d^{d-2}$ such factorisation which coincides with the degree of the $\LL$ map.
\end{rmk}

Using the transitivity of the Hurwitz action and the degree of the $\LL$ map, \cite{DM24} showed the following:

\begin{prop}[Proposition 11.6 \cite{DM24}]\label{prop:MonoAndPol}
	Monic centered polynomials with the same multiset of critical values and the same chains of noncrossing partitions are equal.
\end{prop}
\subsection{The Geometric Combinatorics Map}
Given a poset $\mathcal{P}$, its \textit{order complex}, denoted $\abs{\mathcal{P}}$, is a simplicial complex whose $k$-simplex is a $k$-chain in $\mathcal{P}$. Since a monic centered polynomial is equipped with two chains of noncrossing partitions, \cite{DM24} associates every polynomial to a point in $\abs{\NCP_d}\times \abs{\NCP_d}$, which they called the \textit{geometric combinatorics map}. In this section, we fix $\closedrect=[x_\ell,x_r]\times i[y_b,y_t]$.
\begin{definition}[Barycentric Coordinates]
	Let $x_i$ denote the values of $\re(\cvl(p))$, ordered in increasing order. The $i^{th}$-\textit{relative width} is defined to be $\mathbf{wt}_i=(x_i-x_{i-1})/(x_r-x_\ell)$ where $i=1,2,\cdots,k+1$, $x_0=x_\ell$ and $x_{k+1}=x_r$. The left-right chain of noncrossing partition $[\lambda_1]^L<[\lambda_2]^L<\cdots<[\lambda_{k+1}]^L$ defines a $k$ simplex in the order complex, with vertices labelled by $[\lambda_i]^L$. We define the \textit{barycentric coordinates of $p$ in that simplex to be}: $\textsc{Bary}_{\text{Re}}(p)=\sum_i^{k+1}\mathbf{wt}_i[\lambda_i]^L$. 
	
	The barycentric coordinates for the imaginary part of the critical values is defined analogously: $\textsc{Bary}_\text{Im}(p)=\sum_j^{k+1}\mathbf{wt}_j[\lambda_j]^B$.
\end{definition}

Notice that a small variation on the image point changes the relative widths of the critical values slightly, thus the following is clear. 
\begin{prop}
	The map $\textsc{Bary}_\text{Re}(p):\poly_d^{mc}(\closedrect)\to\abs{\NCP_d}$ is continuous.
\end{prop}

This leads us to one of the main results of \cite{DM24}.
\begin{thm}[Interval Theorem, Theorem A \cite{DM24}]\label{thm:injBaryAn}
	The map $\textsc{Bary}_\text{Re}(p):\poly_d^{mc}(\closedint)\to\abs{\NCP_d}$ is bijective. In particular, the space $\poly_d^{mc}(\closedint)$ is homeomorphic to the order complex $\abs{\NCP_d}$.
\end{thm}

\begin{definition}[Geometric Combinatorics map]
	Let $p$ be a polynomial in $\poly_{d}^{mc}(\closedrect)$, let $\text{Re}(\cvl)$ be the multisets of the real part of the critical values, let $\text{Im}(\cvl)$ be the multisets of the imaginary parts of the critical values, and let $[\lambda_i]^L$, $[\lambda_j]^B$ be the sequence of noncrossing partitions coming from the left and bottom subrectangles respectively. The geometric combinatorics map is defined to be
	\begin{align*}
		\textsc{GeomCom}: \poly_{d}^{mc}(\closedrect)&\to\abs{\NCP_d}\times \abs{\NCP_d}\\
		p&\mapsto (\textsc{Bary}_{\text{Re}}(p),\textsc{Bary}_\text{Im}(p))
	\end{align*}
\end{definition}
\begin{thm}[Rectangular Theorem, Theorem C \cite{DM24}]\label{thm:InjGeomCom}
	The map $\textsc{GeomCom}$ is injective. Thus, $\poly_d^{mc}(\closedrect)$ embeds into $\abs{\NCP_d}\times \abs{\NCP_d}$.
\end{thm}

The left and right sides of $\closedrect$ can be identified using a nonsplitting homotopy. By lifting the homotopy to the space of monic centered polynomials, we obtain the topological space $\poly_d^{mc}(\closedann)$. Let $\mathcal{K}_d$ denote the dual braid complex of the braid group of $d$ strands, \cite{DM24} showed the following:

\begin{thm}[Theorem 15.1 \cite{DM24}]
	The space $\poly_d^{mc}(\closedann)$ embeds into the space $\mathcal{K}_d\times \abs{\NCP_d}$.
\end{thm}

The dual braid complex and the quotient hyperplane complement are two different classifying spaces of the braid group. As polynomials with nonzero critical values are polynomials with distinct roots, the space $\poly_d^{mc}(\closedann)$ is a compactification of the quotient hyperplane complement. \cite{DM24} constructed an explicit homotopy between the two classifying spaces as follows:
\begin{thm}[Theorem B and G, \cite{DM24}]
	The deformation retraction $\closedann$ to $\circcle$ induces a deformation retraction from $\poly_d^{mc}(\closedann)$ to $\poly_d^{mc}(\circcle)$, where the latter is homeomorphic to the dual braid complex. 
\end{thm}
\section{Space of polynomials for type $G(k,1,d)$}
%%%%%%%%%%%%%%%%%%%%%%%%%%%%
In this section, we define the space of polynomials that is homeomorphic to $\C^d/G(k,1,d)$. We also discuss the noncrossing combinatorics of type $B_d$.
\subsection{The groups of type $G(k,1,d)$}
Recall that $G(k,1,d)$ is a finite complex reflection group generated by:
	\begin{align*}
		\begin{pmatrix}
			w_1&0&\cdots&&0\\
			0&w_2&0&\cdots&0\\
			0&0&\ddots&&0\\
			0&&&&w_d
		\end{pmatrix}
	\end{align*}
	where $w_i=\zeta^{k_i}$  for some $k_i\in \{0,1,2,3,\cdots\}$ and $\zeta$ is a primitive $k$-th root of unity, together with all the permutation matrices in $\text{GL}(d,\C)$. A \textit{permutation matrix} is a matrix where every row and every column has only a single entry of $1$ with all other entries $0$. In particular, the group $G(2,1,d)$ is known as the Coxeter group of type $B_d$, which has the following presentation:
	\begin{align*}
		G(2,1,d)\cong \cox(B_d)=\langle s_1,s_2,\cdots,s_d\ \vert\ (s_is_j)^{m_{ij}}=1\rangle	\quad 	m_{ij}=\begin{cases}
			4,& i=1,j=2\\
			3,& i\neq 1, j=i+1\\
			2,&\text{otherwise}
		\end{cases}
		\end{align*}
	\begin{rmk}
		Consider the group $\sym_{2d}$ acting on the set $\{\pm 1,\cdots,\pm d\}$. Then, the group $\cox(B_d)$ can be identified, as a subgroup of $\sym_{2d}$, with the centraliser of the involution sending $i\mapsto -i$. 
	\end{rmk}
	\begin{rmk}
		The invariant ring of finite complex reflection groups are polynomial rings due to Chavalley and Shephard-Todd. Let $e_i(x_1,\cdots,x_d)$ be the $i$-th elementary symmetric polynomial in $d$ variable, it is not difficult to show that the invariant ring of $G(k,1,d)$ is generated by the elements $e_i(x_1^k,\cdots,x_d^k)$.
	\end{rmk}
\subsection{Hyperplane Complement of $G(k,1,d)$}
The hyperplane complement of a finite complex reflection group acting on $\C^d$ is the subspace of $\C^d$ consisting of points with trivial stabiliser. Given a finite complex reflection group, one defines an invariant polynomial $D$ as follows:
\begin{align*}
	D=\prod_{h\in \mathcal{H}}L_h^{o(h)}
\end{align*}
where $\mathcal{H}$ is the set of hyperplanes, $L_h$ is a linear equation which vanishes on the hyperplane $h$, and $o(h)$ is the order of the cyclic group stabilising $h$. Since $D$ is invariant, it is an element in the invariant ring, $\C[x_1,\cdots,x_d]^W$, that is called a \textit{discriminant}, denoted as $\Delta(W;f_1,\cdots,f_d)$, where the $f_i$'s are the invariant polynomials of $W$ ordered in increasing degrees. Under the quotient map $\C^d\to \C^d/W$, the union of the hyperplanes map onto the vanishing locus of $\Delta(W;f_1,\cdots,f_d)$. Thus, one sees that the space $\C^d-\{v\in \C^d\ \vert\ \Delta(v)=0\}$ is homeomorphic to the quotient hyperplane complement. From this, one concludes that if two groups have discriminants which differ by a constant, then the fundamental group of their quotient hyperplane complement must be isomorphic. 
\begin{prop}[Corollary 2.26 \cite{orlik1988discriminants}]\label{prop:isodisc}
	Let $k>2$, then the discriminant of $G(k,1,d)$ differs from that of $\cox(B_d)$ by a constant.
\end{prop}

We see that the problem of studying the quotient hyperplane complement of $G(k,1,d)$ reduces to the study of the same problem for $\cox(B_d)$.

\subsection{Noncrossing partitions of type $B_d$}
The noncrossing partions of $2d$ elements form a lattice denoted as $\NCP_{2d}$. The partitions which are $\pi$-rotationally symmetric form a subposet of $\NCP_{2d}$ which is called the noncrossing partitions of Type $B_d$. This subposet enjoys many of the same properties as $\NCP_d$, for example it is a lattice and is isomorphic to the lattice of reflection factorisation of a Coxeter element in $\cox(B_d)$. A lot of work involving $\NCP(B_d)$ was done by Reiner in \cite{reiner1997non}, where he showed many properties of the lattice which we explain here.

\begin{definition}[Reiner]\label{def:NCPBn}
	The \textit{noncrossing partition of Type $B_d$} is a partition, $[\pi]$, of $\{\pm1,\cdots,\pm d\}$ such that for any block $B$ in $[\pi]$, its negative, $-B$, is also a block of $[\pi]$, and there is at most one block containing both $i$ and $-i$. We denote $\NCP(B_d)$ to be the \textit{poset} of all noncrossing partitions of Type $B_d$, where $[\sigma]<[\tau]$ if each block of $[\sigma]$ is contained in a block of $[\tau]$.
\end{definition}
\begin{example}\label{ex:NCPBn}
	Given the set $\{\pm 1,\cdots\pm n\}$, label $2d$ equally spaced dots on a circle with $d,\ (d-1),\cdots,\ 1,\ -d,\ -(d-1),\cdots,\ -1$ in clockwise order. By labelling the dots in this manner, every noncrossing partition partition of type $B_d$ is equivalent to the partition having a $\pi$-rotational symmetry in their diagrams. As an example, the poset $\NCP(B_3)$ is shown in Figure~\ref{Fig:NCPB3}. 
	\begin{figure}
		\begin{center}
			\includegraphics[width=\textwidth]{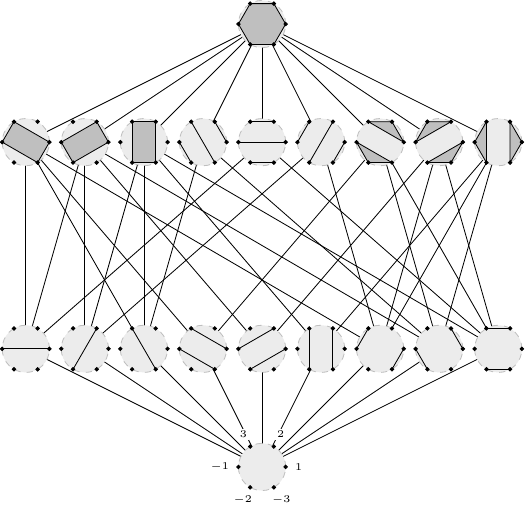}
			\caption{The figure shows $\NCP(B_3)$, the poset of noncrossing partitions of type $B_3$.The dots are labelled by $3,\ 2,\ 1,\ -3,\ -2,\ -1$ in the clockwise direction. Each partition is invariant under a $\pi$ rotation. The lines represent the ordering of the poset, with the top element representing the maximal element of the poset.}\label{Fig:NCPB3}
		\end{center}
	\end{figure}
\end{example}
\begin{rmk}
	\label{rmk:ReinerRotation} An equivalent way to state Definition~\ref{def:NCPBn} is that $\NCP(B_d)$ is the sub-poset of $\NCP_{2d}$, labelled using the convention in Example~\ref{ex:NCPBn}, consisting of partitions which are invariant under a $\pi$ rotation. 
\end{rmk}

The group $\cox(B_d)$ is a reflection group that is isomorphic to the group of symmetries of an $d$-cube $[-1,1]^d$. Being a reflection group, it is generated by a minimal set of reflections whose product is called a Coxeter element. Each Coxeter element has a lattice of reflection factorisation, denoted as $([1,c],\leq_R)$, where for any elements $u,v\in W$,  $u\leq_Rv$ if $\ell_R(v)=\ell_R(u)+\ell_R(u^{-1}v)$. It is well-known that the two lattices are isomorphic.

\begin{prop}[\cite{reiner1997non}]
	$\NCP(B_d)$ is isomorphic to the lattice $([1,c],\leq_R)$. In particular, $\NCP(B_d)$ is a lattice and it contains $d^d$ maximal chains.
\end{prop}

\subsection{Space of Polynomials for $G(k,1,d)$}
\begin{definition}
	Let $k\geq 2$ be an integer. The space $\polyk{k}(\C)$  is the space of polynomials of the form $p(z^k)$, where $p$ is a monic polynomial of degree $d$.
\end{definition}
The following is easy to verify.
\begin{prop}
	$\polyk{k}(\C)$ consists of polynomials of the form:
	\begin{align}
		p(z^k)=z^{dk}+a_{d-1}z^{(d-1)k}+\cdots+a_{1}z^k+a_0
	\end{align}
	Furthermore, if $r_1^k,\cdots,r_d^k$ are the roots of the polynomial $p(z)$, and $\zeta_k^jr_i$ are the roots of $p(z^k)$, then $a_i=e_i(r_1^k,\cdots,r_d^k)$, where $e_i(x_1,\cdots,x_d)$ is the $i^{th}$ elementary symmetric polynomial with $d$ variables.
\end{prop}

\begin{prop}
	$\polyk{k}(\C)$ is homeomorphic to $\C^d/G(k,1,d)$.
\end{prop}
\begin{proof}
	The coefficients of the polynomials in $\polyk{k}(\C)$ are the invariant polynomials of $G(k,1,d)$. It is well-known, see \cite{lehrer2009unitary}, that the invariant polynomials define an orbit map $V\to \C^n$ that descends to a bijection between $V/W\to \C^n$.
\end{proof}

In view of Proposition~\ref{prop:isodisc}, we only have to consider the space $\polyk{2}(\C)$.
\begin{rmk}\label{rmk:FactsOfPoly2}
	We collect some quick facts about $\polyk{2}(\C)$:
	\begin{enumerate}
		\item The roots of any $p\in \polyk{2}(\C)$ have $\pi$ rotation symmetry, that is, if $a$ is a root of $p$, then so is $-a$.
		\item There are $2d-1$ critical points up to multiplicity. The point $z=0$ is always a critical point.
		\item There are $n=d$ critical values up to multiplicity. To see this, notice that $(p(z^2))'=p'(z^2)\cdot 2z$. Thus, if $\alpha_i$ is a critical value, then so is $-\alpha_i$. However, since $p(z^2)$ evaluates to the same value for both $\alpha_i$ and $-\alpha_i$, nonzero critical points must come in pairs that have the same critical value. Since a generic polynomial has $2d-1$ critical points and $0$ is equal to its negative, we have $(2d-2)/2=d-1$ different pairs and $d=n$ critical values. 
	\end{enumerate}
\end{rmk}

\section{The Lyashko-Looijenga Map}
In \cite{DM24}, the authors used a map that relates polynomials to their multiset of critical values. While the same map is also defined on the space $\polyk{2}(\C)$, the resulting map sends polynomials with $d$ coefficients to multisets with $2d-1$ points, which are not the appropriate dimensions. In this section, we find an appropriate generalisation of the map for our space of polynomials, which we denote as  $\LL_B$, which counts the multiplicities of the critical points and critical values in a way that preserves the stratified covering nature of the original map. This map is related to the $\LL$ map defined by Bessis in \cite{bessis2015finite}.

\begin{definition}[Reduced Shape]\label{def:redShape}
	Let $[\lambda]=\{2k_0+1,2k_1,\cdots, 2k_\ell\}$  be an integer partition of $2n-1$. We call the integer partition of $n$, $\{k_0+1,k_1,\cdots,k_\ell\}$ its \textit{reduced integer partition}. Let $p\in\polyk{2}(\C)$, the multiplicities of $\cvl(p)$ can be written as $\{z_0^{2k_0+1},z_1^{2k_1},\cdots,z_\ell^{2k_\ell}\}$, where $z_0=p(0)$. The \textit{reduced shape} of the critical values is $\rshp(\cvl(p))=\{k_0+1,k_1,\cdots,k_\ell\}$.
\end{definition}
\begin{definition}\label{def:LLBn}
	Let $p(z^2)$ be a polynomial in $\polyk{2}(\C)$, and let $\{a_0,\cdots,a_\ell\}$ be its critical values with reduced shape $\{k_0+1,k_1,\cdots,k_\ell\}$. The \textit{Lyashko-Looijenga map of Type }$B_d$, denoted $\LL_B$, is the map sending the polynomial $p(z^2)$ to its multiset of critical values with reduced shape. Concretely:
	\begin{align*}
		\LL_B:\polyk{2}(\C)&\to\multi_{n}(\C)\\
		p(z^2)&\mapsto \{a_0^{k_0+1},a_1^{k_1},\cdots,a_\ell^{k_\ell}\}
	\end{align*}
\end{definition}

\begin{rmk}
	The original $\LL$ map was defined by Lyashko and Looijenga \cite{looijenga1974complement} in the context of singularity theory. The map was then generalised to all well-generated finite complex reflection groups by Bessis in \cite{bessis2015finite} using the discriminant of the groups. Our map here can be viewed as adapting Bessis' $\LL$ map to the context of polynomials and critical values.
\end{rmk}

\begin{rmk}[$n=d$]
	Following \cite{DM24}, the integer $n$ represents the number of critical values, which is equal to $d-1$ in the case of the $\poly_{d}^{mc}(\C)$. However, in the case of $\polyk{m}(\C)$, Remark~\ref{rmk:FactsOfPoly2} and Definition~\ref{def:redShape} tell us that $n=d$.
\end{rmk}

Since the critical values of a polynomial depend only on the coefficients of the polynomial,  $\LL_B$ is an algebraic map. 
\begin{thm}
	The degree of $\LL_B$ is $n^n$.
\end{thm}
The theorem tells us that a generic point in $\multi_n(\C)$ has $n^n$ preimages. The theorem can be proved by emulating the proof found in \cite{lando2004graphs}. This involves applying the B\'{e}zout theorem for quasihomogeneous map, and then showing that the map $\LL_B$ is quasihomogeneous with the appropriate weights.

Recall that Theorem~\ref{thm:LLStratCov} shows that $\LL$ map is a covering map if we restrict the domain and codomain to specified strata. We show that the $\LL_B$ map also has this property. 
\begin{definition}[Reduced shape of critical points]\label{def:CritPtShape}
	Given $p\in\polyk{2}(\C)$, its mutiset of critical points $\cpt(p)$ is always of the form $\{0^{2k_0+1},\alpha_1^{k_1},\cdots,\alpha_\ell^{k_\ell},-\alpha_1^{k_1},\cdots,-\alpha_\ell^{k_\ell}\}$. Notice that the second half of the multiset is technically redundant. Hence, we shall denote $\cpt(p)=\{0^{2k_0+1},\alpha_1^{2k_1},\cdots,\alpha_\ell^{2k_\ell}\}$, with the implicit understanding that the even multiplicities of $\alpha_i$ come from counting their antipodal images. This allows us to define its reduced shape analogously as $\rshp(\cpt(p))=\{k_0+1,k_1,\cdots,k_\ell\}$.
\end{definition}
\begin{definition}[Stratification]
	The space of polynomials has a double stratification by reduced shape of critical points and critical values. Let $\cpt(p)$ and $\cvl(p)$ denote the multiset of critical points and critical values respectively, and let $\lambda=\lambda(p)$ and $\mu=\mu(p)$ be the reduced shape of $\cpt(p)$ and $\cvl(p)$ respectively. We denote $\lambda(p)\xrightarrow{p}\mu(p)$ to represent the polynomial $p$ mapping $\cpt(p)$ with reduced shape $\lambda$ to $\cvl(p)$ with reduced shape $\mu$, this gives us an ordering $\lambda \rightarrow\mu$. Note that we can do this since $\alpha_i$ and $-\alpha_i$ share the same critical value. Let $\poly_{\lambda\to\mu}(\C)=\{p\in\polyk{2}\ \vert\ \lambda\to\mu \text{ is }\lambda(p)\xrightarrow{p}\mu(p)\}$, we get the following stratification of our polynomial space
	\begin{align*}
		\polyk{2}(\C)=\bigsqcup_{\substack{\lambda,\mu\vdash d\\\lambda\rightarrow\mu}}\poly_{\lambda\to\mu}(\C)
	\end{align*}
\end{definition}
\begin{thm}\label{thm:LLStratCovBn}
	Let $\lambda,\mu\vdash n$ such that $\lambda\rightarrow \mu$. The restricted map $\LL_B:\poly_{\lambda\to\mu}^{2}(\C)\to\multi_\mu(\C)$ is a covering map.
\end{thm}
The proof of the theorem is very similar to the proof of Theorem~\ref{thm:LLStratCov}. Since the proof is rather technical, it is included in the appendix. 

In \cite{DM24}, the authors proved that paths in $\multi_n(\C)$ with weakly increasing shape, that is if $\lambda(t)$ is the shape of a path, then for all $s\geq t$, $\lambda(s)\geq \lambda(t)$, have unique path lifiting under the $\LL$ map. As the proof only relies on the stratified covering nature of the map, we have the analogous result.
\begin{thm}\label{thm:PathLiftingLLB}
	Let $\alpha:I\to\multi_n(\C)$ be a path such that $\alpha$ has a weakly increasing reduced shape. Let $\tilde{\alpha}(0)$ be a polynomial in $\polyk{2}(\C)$ such that $\LL_B(\tilde{\alpha}(0))=\alpha(0)$, then there exists a unique lift of $\alpha$ to a path $\tilde{\alpha}:I\to\polyk{2}(\C)$ such that $\LL_B\circ\tilde{\alpha}=\alpha.$ 
\end{thm}
In \cite{DM24}, the authors also introduced the notion of nonsplitting, nonmerging and point-preserving homotopies of $\multi_n(\C)$. In particular, nonsplitting homotopies can only increase the shape of the points inside $\multi_n(\C)$. This allowed them to define a polynomial homotopy by lifting a nonsplitting homotopy of $\multi_n(\C)$. Thus, we have the following:
\begin{definition}[Polynomial homotopy]
	Let $H:U\times I\to \C$ be a nonsplitting homotopy, then there exists a homotopy $\widetilde{H}_n:\polyk{2}(U)\times I\to \polyk{2}(\C)$ satisfying 
	\begin{align*}
		\LL_B\circ \widetilde{H}_n=H_n\circ(\LL_B\times id_I) 
	\end{align*}
	we call $\widetilde{H}_n$ a polynomial homotopy.
\end{definition}
To be precise, $H$ induces a multiset homotopy, and $\widetilde{H}_n$ is defined to be the map such that the following diagram commutes:
\begin{displaymath}
	\begin{tikzcd}
		\polyk{2}(U)\times I\arrow[r,"\widetilde{H}_n"]\arrow[d,"\LL_B\times id_I"]&\polyk{2}(\C)\arrow[d,"\LL_B"]\\
		\multi_{n}(U)\times I\arrow[r,"H_n"]&\multi_{n}(\C)
	\end{tikzcd}
\end{displaymath}
Let $H_n(M,t)=M_t$ be a multiset in $\multi_n(\C)$, then there exists some polynomial $p$ in $\polyk{2}(\C)$ such that $\LL_B(p)=\cvl(p)=M_0$. Furthermore, since $H_n$ is nonsplitting, by the uniqueness of path lifting, we can find a one parameter family of polynomials $p_t$ such that $\LL_B(p_t)=M_t$. We define $\widetilde{H}_n(p,t)=p_t$.

The following results are also from \cite{DM24} with minor changes to their proofs to our definition of $\LL_B$.
\begin{lemma}\label{lem:pointPreservingIsoBn}
	Let $H:U\times I\to \C$ be a point-preserving homotopy and $\widetilde{H}_n$ be the induced polynomial homotopy. Let $t,s\in I$ and $t<s$, then $\widetilde{H}_n$ maps $\polyk{2}(U_t)$ bijectively onto $\polyk{2}(U_s)$.
\end{lemma}
Let $V\subset U$ be subspaces of $\C$ and let $H:U\times I\to \C$ be a nonsplitting point homotopy. $H$ is a \textit{nonsplitting deformation retracting} point homotopy from $U$ to $V$ if the image of $H$ lies in $U$ and $H$, with range restricted to $U$, is a deformation retraction from $U$ to $V$. The following can be adapted from \cite{DM24}.
\begin{prop}[Deformation Retraction]\label{prop:DefRetB}
	If $H$ is a nonsplitting deformation retracting point homotopy from $U$ to $V$ inside $\C$, then there is a $\widetilde{H}_n$-induced deformation retraction from $\polyk{2}(U)$ to $\polyk{2}(V) $ inside $\polyk{2}(\C)$.
\end{prop}

The existence of polynomial homotopies allows one to homotope $\C$ to $\openrect$ via a nonsplitting homotopy, which lifts to a homotopy from $\polyk{2}(\C)$ to $\polyk{2}(\closedrect)$. This gives us the following result. 
\begin{thm}
	$\polyk{2}(\closedrect)$ is the compact closure of $\polyk{2}(\openrect)$.
\end{thm}
\begin{proof}[Sketch of proof]
	The proof is similar to Theorem~\ref{thm:compactification}, where a point-preserving homotopy from $\C$ to the interior of $\closedrect$, denoted as $\openrect$, is constructed, which induces a homotopy from $\polyk{2}(\C)$ to $\polyk{2}(\openrect)$. The $\LL_B$ map can be shown to be a proper map which implies that $\polyk{2}(\closedrect)$ is compact. By constructing a point-preserving homotopy that shrinks $\closedrect$ into $\openrect$, one can show that $\polyk{2}(\closedrect)$ is the smallest closed set containing $\polyk{2}(\openrect)$, and hence must be its compact closure. 
\end{proof}

\section{Noncrossing Combinatorics of a Polynomial}
\subsection{Noncrossing Partitions}
We turn our attention away from the space of polynomials to a single polynomial. Recall that in \cite{DM24}, the authors associate to every polynomial in $\poly_d^{mc}(\closedrect)$ two chains of noncrossing partitions. Since we have the inclusion $\polyk{2}(\closedrect)\subset \poly_{dk}^{mc}(\closedrect)$, the polynomials in $\polyk{2}(\closedrect)$ also have two chains of noncrossing partitions. However, we prove that the noncrossing partitions are restricted to those of type $B_d$.

\begin{definition}[Labels]\label{def:labelsBn}
	Let $R$ be a closed rectangle containing $\closedrect$ in its interior and let $P=p^{-1}(R)$ be the preimage polygon.  Consider the preimage of the right boundary of $R$ together with the positive real ray. There exists a unique preimage such that the preimage of the positive real ray is asymptotic to the real axis, we call that boundary $r_d$ ($r$ for the right side of $R$). Starting from $r_d$, we label the preimages of $r$ by $d,\ (d-1),\cdots,\ 1,\ -d,\ -(d-1),\cdots,\ -1$ in clockwise order. 
	
	We label the other boundary components in similar way. Going clockwise starting from $r_d$, the boundary component adjacent to it, that is the preimage of the bottom component, is labelled $b_d$. The boundary component adjacent to $b_d$, which is the preimage of the left component, is labelled $l_d$ and so on. See Figure~\ref{fig:pullbacklabelBn}.
	
	\begin{figure}\
		\centering
		\includegraphics[width=0.9\linewidth]{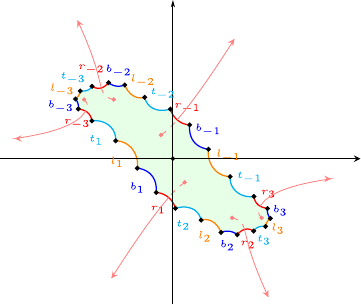}
		\caption{The figure shows the labelling of the sides of the preimage polygon. We first determine the side with the label $r_3$, which is defined to be the side that intersects the preimage of the positive real ray, shown here in light red, which is asymptotic to the real axis in the domain, shown here in black. The other sides are then labelled clockwise in decreasing order until $1$ and then in increasing order from $-3$.}
		\label{fig:pullbacklabelBn}
	\end{figure}
\end{definition}
\begin{rmk}[$2d$ roots of unity]\label{rmk:droots}
	Given any monic polynomial, the positive real axis maps asymptotic to the positive real axis, and, in the other direction, the preimage of the positive real axis consists of $2d$ open half rays asymptotic to the $2d$-roots of unity. Note that if the roots of the polynomial are distinct, the $2d$ open rays start at the $2d$ roots of $p$. 
\end{rmk}

Recall that given a polynomial $p$ with critical values lying inside a fixed $\closedrect$, we choose a closed rectangle $R=I' \times i J'$ such that $\closedrect$ is in its open interior. We subdivide $R$ into a sequence of left subrectangles by subdividing the interval $I'$ such that every subinterval not containing the endpoints contains only one value in $\re(\cvl(p))$. Similarly, we subdivide $R$ into a sequence of bottom subrectangles by subdividing the interval $J'$ so that every subinterval not containing the endpoints has only one value in $\imag(\cvl(p))$. See Example~\ref{ex:AnExample} for how one such subdivision is done. The subrectangles pull back under $p$ to give partitions of the preimages of $l_i$ and $b_j$, the preimages of the left and bottom side of $R$ under $p$ respectively.
\begin{prop}
	The pullback of the subrectangles define a noncrossing partition of type $B_d$.
\end{prop}
\begin{proof}
	Following the work of \cite{DM24}, since the polynomial $p$ is monic and centered, the partition obtained is a noncrossing partition of type $A_{2d}$. By construction, $p$ is the composition of a monic polynomial and the map $z\mapsto z^2$. Since the squaring map sends antipodal points to the same value, the noncrossing partition must have $\pi$ rotational symmetry and hence it is of type $B_d$.
\end{proof}
\begin{cor}
	The sequence of noncrossing partitions of type $B_d$, $[\lambda_i]^L$, obtained from the sequence of left subrectangles $R_1^L\subset R_2^L\subset\cdots\subset R_k^L\subset R$ is strictly increasing from the discrete partition to the indiscrete partition. This is also true for the sequence of noncrossing partitions $[\lambda_i]^B$ obtained from the sequence of bottom subrectangles $R_1^B\subset R_2^B\subset\cdots\subset R_k^B\subset R$.
\end{cor}
\begin{proof}
	By Proposition~\ref{prop:strictInc}, we know that as noncrossing partitions of $2d$ elements, the sequence of noncrossing partitions is strictly increasing from the discrete partition to the indiscrete partition. The previous proposition shows that all the partitions belong in the poset $\NCP(B_d)$.
\end{proof}

\subsection{Monodromy Action} 

Let $\closedrect$ be a fixed closed rectangle, and let $p\in\polyk{2}(\closedrect)$. Let $R\supsetneq \closedrect$ be a closed rectangle containing $\closedrect$ in its interior, and let $P=p^{-1}(R)$. Let $z_a$ denote the bottom left corner of $R$ which we use as a basepoint. The boundary of the right subrectangles $R_i^L$, denoted as $\D R_i^L$, is a based loop with basepoint at $z_a$, which we orient counter-clockwise. Since the boundary curve does not contain any critical values, it can be lifted through $p$ into disjoint paths in $P$. 

Let $\widehat{P}=P-p^{-1}(\cvl(p))$ be the preimage polygon with the preimage of the critical values removed, and let $\widehat{R}=R-\cvl(p)$ be the closed rectangle with the critical values removed. Since $p:\widehat{P}\to\widehat{R}$ is a covering map, based loops with basepoint at $z_a$ can be lifted to paths in $\widehat{P}$ which give a permutation of the preimages of $z_a$. This gives us a group homomorphism known as the monodromy action:
\begin{align*}
	\textsc{Mono}: \pi_1(\widehat{R},z_a)\to \sym_{2n}
\end{align*}
The image of $\D R_i^L$ in $\sym_{2n}$ can be obtained by replacing the $k$-blocks in the nonrossing partition $[\lambda_i]^L$ by $k$-cycles, where elements of the $k$-cycle are listed in counter-clockwise order according to the diagram of $[\lambda_i]^L$. Given the rotational symmetry of the noncrossing partitions, the following is clear. Given the rotational symmetry of the noncrossing partitions, the following is clear.
\begin{prop}
	The image of $ \pi_1(\widehat{R},z_a)$ under $\textsc{Mono}$ is contained in the centraliser of the involution $\iota(i)=-i$. In particular, the image is contained in a subgroup of $\sym_{2d}$ isomorphic to $\cox(B_d)$.
\end{prop}

The maximum element in the chains of noncrossing partitions gives a $2d$-cycle, which is a Coxeter element of $\cox(B_d)$. Let $\pi_i^L=\textsc{Mono}(\D R_i^L)$, the sequence of left subrectangles gives us a chain of factorisation of the Coxeter element, $\pi_1^L<\pi_2^L<\cdots<\pi_k^L$. Given such a chain, we define the elements $\sigma_i=(\pi_{i-1}^L)^{-1}\pi_{i}^L$ for $i=2,\dots,k$, and $\sigma_1=\pi^L_1$. These elements give us a (left-right) factorisation of the $2d$ cycle, denoted as $[ \sigma_1\ \cdots\ \sigma_{k-1}]$. Similarly, if we let $\pi_j^B=\textsc{Mono}(\D R_j^B)$, the elements $\tau_j=(\pi_{j-1}^B)^{-1}\pi_{j}^B$ for $j=2,\dots,\ell$ and $\tau_1=\pi_1^B$ define a (bottom-top) factorisation, $[ \tau_1\ \cdots\ \tau_{\ell-1}]$, of the same $2d$ cycle. Since all the $\pi_i^L$ and $\pi_j^B$ are in $\cox(B_d)$, the following is clear. 
\begin{prop}
	The sequences $[ \sigma_1\ \cdots\ \sigma_{k-1}]$ and $[ \tau_1\ \cdots\ \tau_{\ell-1}]$ are factorisations of the $2d$ cycle, a Coxeter element in $\cox(B_d)$. In other words, let $c$ denote the Coxeter element, then $c=\sigma_1\ \cdots\ \sigma_{k-1}=\tau_1\ \cdots\ \tau_{\ell-1}.$
\end{prop}
\begin{rmk}[Reflection Factorisation]\label{rmk:refFactBn}
	When the reduced shape of the critical values is $1^n$, the factorisation of the Coxeter element is a reflection factorisation.
	To see this, notice that the multiplicity of a critical value is equal to the reflection length, as an element in $\sym_{2d}$, of the monodromy action around that critical value. One can also check that the reflection length of the generators of $\cox(B_d)$, when embedded inside $\sym_{2d}$, have shape $2^{n-1}1$. Therefore, as the reduced shape $1^n$ tells us that the actual shape of the critical values is $2^{n-1}1$, the factorisation of the $2d$-cycle inside $\cox(B_n)$ must have maximal length, i.e. it is a reflection factorisation. 
\end{rmk}
In order to prove that polynomials in $\polyk{2}(\closedrect)$ are determined by their critical values and their monodromy actions, we can emulate the proof found in \cite{DM24} which involves the degree of $\LL_B$ and Hurwitz transitivity on the reflection factorisation of the Coxeter element.
\begin{lemma}
	Let $p,p'\in\polyk{2}(\closedrect)$ such that $\LL_B(p)=\LL_B(p')$, if their monodromy actions coincide, then their left-right and bottom-top factorisation of $c$ are equal.
\end{lemma}
\begin{proof}
	By perturbing the critical values such that their monodromy action remains unchanged, we can assume that they have distinct real parts. The monodromy of a polynomial is then uniquely determined once we choose a basis of the free group. One can choose a basis $\gamma_i$ such that it is homotopic to the path $\D (R^L_{i-1})^{-1}\circ \D R^L_i$ if $i>1$ and $\D R^L_i$ if $i=1$. Since $p$ and $p'$ have the same set of critical values, and their monodromy actions are the same, the monodromy action of $\D R_i^L$ must coincide for all subrectangles as well. In particular, this shows that their left-right factorisations of $c$ must be equal. The same argument works for the bottom-top factorisation. 
\end{proof}
\begin{thm}\label{thm:MonodromyBn}
	Let $p,p'\in\polyk{2}(\closedrect)$, if $\LL_B(p)=\LL_B(p')$ and the monodromy action of $p$ and $p'$ are equal, then $p=p'$.
\end{thm}
\begin{proof}[Sketch of proof]
	Since the set of critical values of $p$ and $p'$ are equal, and their monodromy actions coincide, the previous lemma tells us that their left-right and bottom-top factorisations of the Coxeter element must be the same. 
	
	In \cite{DM24}, the authors constructed a pair of cell complexes known as the regular point complex and the regular value complex. The regular value complex is a cell complex on $R$, where the cells are obtained by superimposing all boundary of the subrectangles onto $R$. Its preimage gives the regular point complex, which is a cell complex on the preimage polygon $P$. Importantly, the cells in the regular point complex which are not `rectangles' must contain a critical point, whose multiplicity is given by the number of sides the cell contains divided by 4. Moreover, the multiplicities of each critical value can also be determined by summing the multiplicities of the critical points.  As their factorisations are the same, their regular point complex must be the same. This tells us that the multiplicities of their critical values must be the same, i.e. $\LL(p)=\LL(p')$. By Proposition~\ref{prop:MonoAndPol}, $p=p'$. 
\end{proof}
So far, we looked at the monodromy action of polynomials in $\polyk{2}(\closedrect)$ by either viewing the group as living inside $\sym_{2d}$ or the polynomial as living inside $\poly_{2d}^{mc}(\closedrect)$. Naturally, one might wonder if there are other polynomials in $\poly_{2d}^{mc}(\closedrect)$ that has monodromy action living inside $\cox(B_d)$. The rotational symmetry of the noncrossing partitions turns out to characterise the polynomials in $\polyk{2}(\closedrect)$.
\begin{lemma}\label{lem:lengthReduc}
	Let $[g_1\ \cdots\ g_k]$ be any factorisation of the Coxeter element in $\cox(B_d)$. When $g_i$'s are viewed as elements in $\sym_{2d}$, their reflection lengths form a multiset which must be reducible, i.e. all but one of the reflection lengths are even.
\end{lemma}
\begin{proof}
	There are two kinds of reflections in $\cox(B_d)$: the first kind are reflections that send $i\mapsto -i$, and the second kind are reflections that send $i\mapsto j$, where $j\neq-i$. Since $\cox(B_d)$ is the centraliser of the involution $\iota$, reflections of the first kind must be of the form $(-i,i)$ while the second kind must be of the form $(ij)(-i,-j)$. In particular, we see that the reflections have lengths equal to $1$ or $2$. 
	
	The Coxeter element is defined to be the product of all the simple reflections. One can check that amongst all the simple reflections, all but one have length $2$. Therefore, the multiset of reflection lengths of the simple reflections must be reducible. Moreover, since conjugation does not change the cycle type, and the Hurwitz action is transitive on the set of all reflection factorisations, we can also conclude that the multiset of reflection lengths of any reflection factorisation must be reducible.
	
	Finally, for non-reflection factorisation of the Coxeter element, the elements $g_i$'s must be a product of reflections. Since the reflection lengths add when one multiplies, the resulting multiset of reflection lengths must also be reducible.
\end{proof}
\begin{prop}\label{prop:piSymBn}
	Let $p\in\poly_{2d}^{mc}(\closedrect)$ be any monic centered polynomial of degree $2d$, and let $[\lambda_1]^L<\cdots<[\lambda_k]^L$ and $[\lambda_1]^R<\cdots<[\lambda_\ell]^R$ be its two chains of noncrossing partitions. Suppose $k,\ell\leq d$, and all the partitions $[\lambda_i]^R$ and $[\lambda_j]^B$ are $\pi$-rotationally symmetric, then $p$ belongs in the subspace $\polyk{2}(\closedrect)$.
\end{prop}
\begin{proof}
	Let $p(z)$ be a polynomial with $\pi$-rotational symmetric noncrossing partitions. Consider the polynomial $p(-z)$ whose domain is obtained from $p(z)$ by a $\pi$-rotation. Since the polynomial has degree $2d$, the preimage of the positive real axis has $2d$ rays asymptotic to the $2d$-roots of unity, see Remark~\ref{rmk:droots}. In particular, the distinguished side labelled $r_d$ and the side labelled $-r_d$ intersect rays which are asymptotic to the positive and negative real axis respectively. Thus, the noncrossing partitions of $p(-z)$ can be obtained from $p(z)$ via a $\pi$-rotation, so they must have the same combinatorics. Finally, since $p(-z)$ has the same multiset of critical values as $p(z)$, Proposition~\ref{prop:MonoAndPol} tells us that $p(z)=p(-z)$.
	
	Since $p(z)$ has even degree, the polynomial $q(z)=p(z)-p(-z)=0$ consists only of odd degree terms. As $q(z)$ is identically zero everywhere, the coefficients of $q$, $2a_{2i+1}$ where $a_{2i+1}$ are the coefficients of $p(z)$, must be zero. Therefore, $p(z)\in\polyk{2}(\closedrect)$.
\end{proof}

\section{Cell Structure on $\polyk{2}(\closedrect)$}
In \cite{DM24}, see also \cite{DM26Mult}, the authors defined a bisimplicial cell structure on $\multi_n(\closedrect)$ that can be lifted via the $\LL$ map to give a cell structure on $\poly_{d}^{mc}(\closedrect)$. The bisimplicial cell structure on $\multi_{n}(\closedrect)$ is defined by first defining the cell structure of $I^n\times J^n$ to be the pull back of the bisimplices in $\multi_n(I)\times \multi_n(J)$, followed by taking the quotient of the bisimplicial cell structure of $I^n\times J^n$ by the diagonal symmetric group action. 
\begin{displaymath}
	\begin{tikzcd}
		\closedrect^n=I^n\times J^n\arrow[r]\arrow[dr]	& \multi_n(I)\times \multi_n(J)\\
		&\multi_n(\closedrect)\arrow[u]
	\end{tikzcd}
\end{displaymath}
The bisimplicial cell structure on $\multi_n(\closedrect)$ is stratified, where points in the same open cell have the same shape. Using the fact that $\LL_B$ is a stratified covering map, these cells can be lifted, stratum by stratum, to give a bisimplicial cell structure on $\polyk{2}(\closedrect)$. Note that since $\multi_n(\closedrect)$ has $n!$ number of top-dimensional simplices, the degree of $\LL_B$ tells us that $\polyk{2}(\closedrect)$ has $d^dn!=n^n\cdot n!$ number of top-dimensional simplices. Finally, by construction of the cells, $\LL_B$ is a cellular map.
\section{Geometric Combinatorics Map}
Recall that in \cite{DM24}, the authors defined the geometric combinatorics map by giving a polynomial barycentric coordinates in $\abs{\NCP_d}\times \abs{\NCP_d}$ obtained by the relative widths of the critical values as well as their chains of noncrossing partitions. Their construction extends to our space of polynomials, defining barycentric coordinates for polynomials inside $\abs{\NCP(B_d)}\times\abs{\NCP(B_d)}$, where $\abs{\NCP(B_d)}$ is the order complex of the lattice of noncrossing partitions of type $B_d$.
\begin{definition}[Barycentric Coordinates for type $B_d$]
	Let $x_i$ denote the values of $\re(\cvl(p))$, ordered in increasing order. The $i^{th}$-\textit{relative width} is defined to be $\mathbf{wt}_i=(x_i-x_{i-1})/(x_r-x_\ell)$ where $i=1,2,\cdots,k+1$, $x_0=x_\ell$ and $x_{k+1}=x_r$. The chain of noncrossing partition $[\lambda_1]^L<[\lambda_2]^L<\cdots<[\lambda_{k+1}]^L$ defines a $k$ simplex in the order complex $\abs{\NCP(B_d)}$, with vertices labelled by $[\lambda_i]^L$. We define the \textit{barycentric coordinates of $p$ in that simplex to be}: $\textsc{Bary}_{\text{B;Re}}(p)=\sum_i^{k+1}\mathbf{wt}_i[\lambda_i]^L$. 
	
	We define the barycentric coordinates for the imaginary part of the critical values analogously. We denote the map as $\textsc{Bary}_\text{B;Im}(p)=\sum_j^{k+1}\mathbf{wt}_j[\lambda_j]^B$.
\end{definition}

\begin{lemma}\label{lem:CellularBn}
	Let $\closedint$ be a closed interval, the map $\textsc{Bary}_{B;\text{Re}}:\polyk{2}(\closedint)\to\abs{\NCP(B_d)}$ is cellular. 
\end{lemma}
\begin{proof}
	We first note that since the barycentric coordinates record all the relative widths of the critical values, the image of $\textsc{Bary}_{B;\text{Re}}$ allows us to recover the multiset of critical values in $\multi_{n}(\closedint)$. Moreover, points in an open cell of $\abs{\NCP(B_d)}$ have the same number of nonzero relative widths which tells us that their critical values have the same shape. Therefore, the map $\abs{\NCP(B_d)}\to\multi_n(\closedint)$ sends open cells to open cells of the same dimension. Finally consider the commutative diagram below:
	\begin{displaymath}
		\begin{tikzcd}
			\polyk{2}(\closedint)\arrow[r,"\textsc{Bary}_{B;\text{Re}}"]\arrow[rd,swap,"\LL_B"]&[2em] \abs{\NCP(B_d)}\arrow[d]\\
			&\multi_n(\closedint)
		\end{tikzcd}
	\end{displaymath}
	Since the cells in $\polyk{2}(\closedint)$ are defined by lifting cells in $\multi_{n}(\closedint)$ through the $\LL_B$ map, and open cells in $\abs{\NCP(B_d)}$ are mapped to open cells in $\multi_{n}(\closedint)$, the map $\textsc{Bary}_{B;\text{Re}}$ is cellular.
\end{proof}
\begin{thm}[Interval Theorem]\label{thm:intervalcase}
	The map $\textsc{Bary}_{B;\text{Re}}:\polyk{2}(\closedint)\to\abs{\NCP(B_d)}$ is bijective. Hence, $\polyk{2}(\closedint)$ is homeomorphic to $\abs{\NCP(B_n)}$.
\end{thm}
\begin{proof}
	As any continuous bijection between a compact space and a Hausdorff space is a homeomorphism, it suffices to show that $\textsc{Bary}_{B;\text{Re}}$ is bijective. Let us first show surjectivity. By using the barycentric coordinates of the top-dimensional simplex, a point in $\abs{\NCP(B_d)}$ defines a multiset of points in $\closedint$. Moreover, the point also determines a chain of noncrossing partitions in $\NCP(B_d)$ and hence a chain of noncrossing partitions in $\NCP_{2d}$. Using Theorem~\ref{thm:injBaryAn}, there exists a polynomial $p\in\poly_{2d}^{mc}(\closedint)$ having the appropriate critical values as well as chain of noncrossing partitions. However, Proposition~\ref{prop:piSymBn} tells us that this polynomial must belong in the subspace $\polyk{2}(\closedint)$.
	
	Suppose now we have $p,q\in\polyk{2}(\closedint)$ such that $\textsc{Bary}_{B;\text{Re}}(p)=\textsc{Bary}_{B;\text{Re}}(q)$. On one hand, their \textit{set} of critical values must coincide. On the other hand, as $p$ and $q$ are sent to the same simplex, their noncrossing combinatorics must be equal. Moreover, their noncrossing combinatorics, together with the relative widths, allow us to work out the reduced multiplicities of the critical values. Since everything must coincide, by Theorem~\ref{thm:MonodromyBn}, $p=q$.
\end{proof}
%	On the other hand, by viewing their noncrossing partitions as living inside $\NCP_{2d}$, Lemma~\ref{lem:CriPointDetCriVal} (\ben{Check if lemma mentions multiplicity (It doesn't right now but I will fix it)}) allows us to recover the multiplicities of the critical values. However, since their noncrossing combinatorics are the same, their multiplicities, and consequently their \textit{multiset} of critical values, are equal. Finally, by Theorem~\ref{thm:MonodromyBn}, $p=q$.
\begin{definition}[Geometric Combinatorics map]
	Let $p$ be a polynomial in $\polyk{2}(\closedrect)$, the geometric combinatorics map of type $B$ is defined to be
	\begin{align*}
		\textsc{GeomCom}_B: \polyk{2}(\closedrect)&\to\abs{\NCP(B_d)}\times \abs{\NCP(B_d)}\\
		p&\mapsto (\textsc{Bary}_{B;\text{Re}}(p),\textsc{Bary}_\text{B;Im}(p))
	\end{align*}
\end{definition}
\begin{thm}[Rectangular Theorem]\label{thm:rectcase}
	The map $\textsc{GeomCom}_B$ is injective. Hence, $\polyk{2}(\closedrect)$ is homeomorphic to its image under $\textsc{GeomCom}$.
\end{thm}
\begin{proof}[Sketch of proof]
	Suppose $p$ and $q$ are two polynomials in $\polyk{2}(\closedrect)$ such that their images under $\textsc{GeomCom}_B$ are equal, then their \textit{set} of critical values and chains of noncrossing partitions must be equal. Recall that the noncrossing partitions build a regular point complex, where preimages of the boundaries of all the subrectangles subdivide $P$ into cells. Moreover, the regular point complex completely determines the multiplicities and positions of the critical values. As the multiplicities of their critical values are all equal, $\LL(p)=\LL(q)$ and by Proposition~\ref{prop:MonoAndPol}, $p=q$.
	%	Recall that the regular value complex is the subdivided rectangle by the lines used to define the left and bottom subrectangles; the regular point complex is the preimage of the regular value complex, in particular it subdivides the $8n$-gons into $4k$-gons where $k>1$. Lemma~\ref{lem:CriPointDetCriVal} tells us that regular point complex determines the position of the critical values. 
	%	
	%	Let $p,q\in\poly_{n}^2(\closedrect)$ such that $\textsc{GeomCom}_B(p)=\textsc{GeomCom}_B(q)$. By definition, $p$ and $q$ have the same left-right and bottom-top noncrossing combinatorics as well as the relative widths. Therefore, the regular point complex of $p$ and $q$ must coincide. This determines the position of all the critical values as well as their multiplicities. By Proposition~\ref{prop:MonodromyBn}, $p=q$.
\end{proof}
\section{Basketball Complex}
In \cite{DM24}, the authors showed that the image of monic centered polynomials under the geometric combinatorics map is the full subcomplex of \textit{basketball vertices}. In this section, we explain the basketball vertices, and we show that the image of $\polyk{2}(\closedrect)$ is also the full subcomplex of basketball vertices in $\abs{\NCP(B_d)}\times \abs{\NCP(B_d)}$.

\begin{definition}[Noncrossing Matching]
	A top-bottom noncrossing matching is set of noncrossing arcs where each arc has an endpoint in $t_i$ and an endpoint in $b_j$ and every $t_i$ is connected by one such arc to exactly one $b_j$. A left-right noncrossing matching is defined analogously.
\end{definition}

Every vertex of $\abs{\NCP_d}\times \abs{\NCP_d}$ is defined by a noncrossing partition of $l_i$'s and a noncrossing partition of $b_i$'s. A noncrossing partition of $l_i$'s is defined by the pullback of a left subrectangle whose boundary line connecting $t$ and $b$ must lift to $2d$ disjoint noncrossing arcs connecting $t_i$'s and $b_i$'s. Thus, a left-right noncrossing partition gives us a top-bottom noncrossing matching. The bottom-top noncrossing partitions behave similarly.

The image of $\textsc{GeomCom}_B$ turns out to be points whose left-right noncrossing matchings and bottom-top noncrossing matchings form `basketballs' when superimposed on one another. 

\begin{definition}[Basketball Vertices]\label{def:bballVert}
	A vertex of $\abs{\NCP_d}\times \abs{\NCP_d}$ is called \textit{basketball} if when the top-bottom noncrossing matching is superimposed on the left-right noncrossing matching, every arc in the top-bottom matching intersects with exactly one arc in the left-right matching. The term basketball came from \cite{martin2007harmonic}.
\end{definition}
 As $\polyk{2}(\closedrect)$ is a subspace of $\poly_{2d}^{mc}(\closedrect)$, and $\NCP(B_d)$ is a sublattice of $\NCP_{2d}$, it is clear that the image of $\polyk{2}(\closedrect)$ is contained in some kind of basketball complex.  This basketball complex turns out to be precisely the full subcomplex consisting of basketball vertices in $\abs{\NCP(B_d)}\times \abs{\NCP(B_d)}$.
\begin{prop}
	The image of $\polyk{2}(\closedrect)$ under the map $\textsc{GeomCom}_B$ is the full subcomplex of $\abs{\NCP(B_d)}\times\abs{\NCP(B_d)}$ consisting of basketball vertices. 
\end{prop}
\begin{proof}
	Since $\polyk{2}(\closedrect)\subset\poly_{2d}^{mc}(\closedrect)$ and $\NCP(B_d)\subset \NCP_{2d}$, it is clear that every pair of noncrossing matching must be basketball by work of \cite{DM24}. What we have to show is that every pair of rotationally symmetric basketball matching comes from a polynomial in $\polyk{2}(\closedrect)$. 
	
	A point in $\abs{\NCP(B_d)}\times\abs{\NCP(B_d)}$ belongs in the basketball complex if any choice of left-right matching and bottom-top matching from the two chains of noncrossing partitions form a basketball. Given such a point, by viewing the noncrossing partitions as living in $\NCP_{2d}$, we obtain a polynomial in $\poly_{2d}^{mc}(\closedrect)$ that has the prescribed combinatorics. However, since the noncrossing combinatorics chosen were $\pi$ rotationally symmetric, by Proposition~\ref{prop:piSymBn}, the polynomial must lie inside $\polyk{2}(\closedrect)$.
\end{proof}
%\begin{rmk}[Metric Cell Complex]
%	In \cite{DM24}, they actually showed that $\poly_{d}^{mc}(\closedrect)$ isometrically embeds into $\abs{\NCP_n}\times \abs{\NCP_d}$, with the latter given a orthoscheme metric defined by CITE.
%\end{rmk}
%%%%%%%%%%%%%%%%%%%%%%%%%%
\section{Branched Annular Complex}
%%%%%%%%%%%%%%%%%%%%%%%%%%
We are now ready to relate the quotient compactified hyperplane complement of $\cox(B_d)$ with its dual Artin complex. In \cite{DM24}, the quotient compactified hyperplane complement of the symmetric group is homeomorphic to the space $\poly_{d}^{mc}(\closedann)$, which embeds into the space $\abs{\NCP_d}\times\mathcal{K}_d$ via $\textsc{GeoComb}$. We will show that this extends to the case of $\cox(B_d)$ as well.
\begin{prop}
	The space $\polyk{2}(\closedann)$ is a quotient of $\polyk{2}(\closedrect)$.
\end{prop}
\begin{proof}[Sketch of proof]
	The proof here is exactly the same as the one in \cite{DM24}, where a point preserving homotopy is applied to $\closedrect$ to identify the left and right sides to give $\closedann$. Meanwhile, the homotopy lifts to a homotopy in $\polyk{2}(\closedrect)$, and polynomials with some critical values lying on the left and right lines are eventually identified.
\end{proof}
The dual Artin complex of the Artin group $\art(B_d)$ is defined to be the quotient of the order complex $\abs{\NCP(B_d)}$ where two simplices with the same oriented edge labels are identified. In addition, every simplex in $\mathcal{K}_{B_d}$ has a standard representative where its edge labels are $[\pi^L\ \sigma_1\ \cdots \sigma_k\ 1]$.

%\begin{thm}
%	$\poly_{n}^2(\closedann)$ is the quotient of $\poly_{n}^2(\closedrect)$.
%\end{thm}
%\begin{prop}
%	$\poly_{n}^{2}(\closedann)$ is a quotient of $\poly_{n}^{2}(\closedrect)$.
%\end{prop}
%\begin{proof}
%	The proof of this is essentially the same as in Proposition~\ref{prop:quotientAnnulus}, since all it relied on was the homotopy lifting of $LL$ map, which is true for $LL_B$.
%\end{proof}

\begin{thm}[Annular Theorem]
	$\polyk{2}(\closedann) $ is homeomorphic to $\polyk{2}(\closedrect)$ with face identification. In particular, $\polyk{2}(\closedann)$ is embeds into the space $\mathcal{K}_{B_d}\times\abs{\NCP(B_d)} $.
\end{thm}
\begin{proof}[Sketch of proof]
		See \cite{DM24} Theorem 15.1 for the full proof. The main idea is to carefully compute the factorisations of the polynomials when the left and right sides are identified, which shows that the polynomials that are identified have representatives that have left-right factorisation exactly the same as the standard representatives of the dual Artin complex. 
	%	The proof of this can be found in Theorem 13.3 and Theorem 15.1 of \cite{DM24}, we shall give a sketch of the argument here. As argued in Proposition~\ref{prop:quotientAnnulus}, only polynomials with critical values on the left and right edge of the rectangle are identified. 
	%	
	%	Consider a simpler case where the critical values lie along a line, $[x_\ell,x_r]$, and there is a homotopy bringing the line to a circle. Suppose we have a polynomial with critical values lying on $x_\ell$ and $x_r$ giving us a factorisation $[\pi^L\ \sigma_1\ \cdots\ \sigma_k\ \pi^R]$. When $x_\ell$ gets identified with $x_r$, the number of distinct critical values drop by one. The new factorisation becomes $[(\delta\pi^R\delta^{-1})\pi^L\ \sigma_1\ \cdots\ \sigma_k\ 1]$. The loop $(\delta\pi^R\delta^{-1})\pi^L$ is homotopic to a clockwise loop containing their merged critical points. Since the identification used is the same as that used to create the dual braid complex, we see that polynomials with critical values lying on a circle is obtained by face identification of the space with critical values lying on $[x_\ell,x_r]$. For the rectangular case, the argument is similar. 
\end{proof}
\begin{rmk}
	In \cite{DM24}, the authors conjectured (Conjecture 3) that their construction would generalise to other types of Coxeter groups.  In particular, they conjectured that given some Coxeter group $W$, there exists a space built out of algebraic factorisation of the Coxeter element, which they called $\textsc{Br}^m_{W,\delta}(\closedann)$, is a compactification of the quotient hyperplane complement of $W$. Our work here showed that this is true for the class of finite complex reflection groups $G(k,1,d)$, where $\polyk{2}(\closedann)$ is a compactification of the quotient hyperplane complement of $G(k,1,d)$.
	\end{rmk}

 Let $\closedann\to \circcle$ be some nonsplitting deformation retraction of the closed annulus to a circle. By Proposition~\ref{prop:DefRetB}, the deformation retract gives a deformation retraction from $\polyk{2}(\closedann)$ to $\polyk{2}(\circcle)$. 
\begin{cor}
	The deformation retraction of $\polyk{2}(\closedann)$ to $\polyk{2}(\circcle)$ gives a deformation retraction of the image of $\polyk{2}(\closedann)$ in $\mathcal{K}_{B_d}\times \abs{\NCP(B_d)}$ to $\mathcal{K}_{B_d}$. Thus, giving a homotopy from the quotient hyperplane complement to the dual Artin complex of type $B_d$.
\end{cor}
\begin{proof}
	It is clear that the deformation retract only affects the bottom-top factorisation of the Coxeter element. In particular, it collapse all factorisations into $1<c$ and all barycentric coordinates to a single point. By composing the deformation retraction with $\textsc{GeomCom}_B$, we can a deformation retraction from the image of $\polyk{2}(\closedann)$ to $\mathcal{K}_{B_d}\times\{\text{point}\}$, which is homeomorphic to the dual Artin complex.
\end{proof}

\printbibliography
\newpage
\appendix
\section{Stratified Covering}
To show that $\LL_B$ is a stratified covering map, we first introduce the following integral.
\begin{align}\label{eqn:CriPointInt}
	p(z)=\int_0^zw^{2k_0+1}(w^2-z_1^2)^{k_1}\cdots (w^2-z_m^2)^{k_m}\ dw
\end{align}
The above integral takes a multiset of points $M=\{0^{2k_0+1},z_1^{2k_1},\cdots,z_m^{2k_m}\}$ with reduced shape $\rshp(M)=\{k_0+1,\cdots,k_m\}$ to a polynomial in $\polyk{2}(\C)$ with those points as critical points counted with the appropriate multiplicites. Following \cite{dougherty2020critical}, we define a map sending the nonzero critical points $\{z_1,\cdots,z_m\}$ to their corresponding critical values:
\begin{align*}
	\theta:\C^m&\to\C^m\\
	(z_1,\cdots,z_m)&\mapsto (p(z_1),\cdots,p(z_m))
\end{align*}
The following is shown in Appendix~ \ref{appendix:crit}
\begin{thm}
	The determinant of the Jacobian of $\theta$ factors as:
	\begin{align}\label{eqn:JacFac}
		\det J(\theta)=C\left(\prod_{j\in [m]}(z_j)^{2(k_0+k_j)+1}\right)\left(\prod_{\substack{i,j\in[m]\\ j\neq k}}(z_k-z_j)^{2k_j}\right)
	\end{align} 
	where $C$ is a nonzero constant. 
\end{thm}
The above determinant identity tells us that $\theta$ is a local homeomorphism when the critical points are distinct. 
\begin{rmk}\label{rmk:uniqueRepBn}
	Any $p\in\polyk{2}(\C)$ can be uniquely represented as an integral of the form \ref{eqn:CriPointInt}. This is because if the multiset of critical points of $p$ is given by $\cpt(p)=\{0^{2k_0+1},z_1^{2k_1},\cdots,z_m^{2k_m}\}$, then we have the equality:
	\begin{align}\label{eqn:basePair}
		p(z)=d\int_0^zw^{2k_0+1}(w^2-z_1^2)^{k_1}\cdots (w^2-z_m^2)^{k_m}\ dw+c
	\end{align}
	with $d\in \C$ a nonzero constant such that the integral returns a monic polynomial, and $p(0)=c$. 
\end{rmk}
Given any $[\lambda']\vdash [2d-1]$, it might not follow that $[\lambda']$ is a shape of $\cpt(p)$ for some $p\in\polyk{2}(\C)$. In particular, for $[\lambda']$ to be a shape of $\cpt(p)$, it must be `reducible', i.e. it must have one odd integer with all other entries even. Given such a $[\lambda']$, the odd integer must correspond to the multiplicity of the critical point $0$, with the other integers being the multiplicities of the other nonzero critical points. 
\begin{definition}
	Given a partition $[\lambda']\vdash[2d-1]$, we can defined two subspaces of $\C^{2d-1}$ by the shape of the points. More specifically, for any $v\in \C^{2d-1}$, we can view $v$ as a multiset, $[v]$, by taking the quotient by the $\sym_{2d-1}$ action on the coordinates. This multiset has a shape which we shall denote $\textsc{Shape}(v)=\textsc{Shape}([v])$. The two subspaces are:
	\begin{align*}
		\C_{[\lambda']}&=\{v\in\C^{2d-1}\ \vert\ \textsc{Shape}(v)=[\lambda']\}\\
		\C^{[\lambda']}&=\{v\in\C^{2d-1}\ \vert\ \textsc{Shape}(v)\geq[\lambda']\}
	\end{align*}
\end{definition}
\begin{prop}\label{prop:critPointstoPolyBn}
	Let $[\lambda']\vdash[2d-1]$, be such that $[\lambda']$ contains exactly one odd integer. Denote $\poly_{[\lambda]}(\C)$ to be the subspace of polynomials having critical points with reduced shape $[\lambda]$ obtained by reducing $[\lambda']$, then there exists a map $\C_{[\lambda']}\to\poly_{[\lambda]}(\C)$ obtained by sending $v\in \C_{[\lambda']}$ to the polynomial of the form \ref{eqn:basePair}, where $c$ is the unique element in $v$ that has odd multiplicity. Moreover, this map is a covering map.
\end{prop}
\begin{proof}
	The above is obvious given the prior discussion except for the claim that it is a covering map. Let $k$ be the number of blocks in $[\lambda']$, then $\C_{[\lambda']}$ can be identified with open subset, $U$, of $\C^k$. To be precise, if the multiset $[v]=\{c^{2k_0+1},z_1^{2\ell_1},\cdots,z_{k-1}^{2\ell_{k-1}}\}$, then we can identify $v$ to the tuple $(c,z_1,\cdots,z_{k-1})$. The map sending points in $U$ to $\poly_{[\lambda]}$ is a map sending $c$ to $c$ while applying the $\sym_{k-1}$ action on the last $k-1$ entries, i.e.
	\begin{align*}
		\varphi:U&\to \C\times \multi_{k-1}(\C)\\
		(c,z_1,\cdots,z_{k-1})&\mapsto (c,[(z_1,\cdots,z_{k-1})])
	\end{align*}
	Since the entries are assumed to be distinct, the map above is a covering map. Finally, $\C\times \multi_{k-1}(\C)$ is homeomorphic to $\poly_{[\lambda]}(\C)$ by Equation~\ref{eqn:basePair}.
\end{proof}

We define the augmented map:
\begin{align*}
	\Theta:\C^{m+1}&\to\C^{m+1}\\
	(c,z_1,\cdots,z_m)&\mapsto (c,p(z_1),\cdots,p(z_m))
\end{align*} 
where $p(z)$ is the integral defined in Equation~\ref{eqn:basePair}. This map sends the set of critical points to their critical values, while also assigning the value to $p(0)$. If we take the partial derivative of the first map ($c\mapsto c$) with respect to the coordinates $\{c,z_1,\cdots,z_m\}$, we get a row consisting of $1$ in the first position and $0$'s everywhere else. Therefore, the Jacobian of $\Theta$ with respect to the coordinates $\{c,z_1,\cdots,z_m\}$ factorises as in \ref{eqn:JacFac} up to a constant.

%\begin{definition}
%	We denote $\poly_{n}^{2,0\to c}(\C)=\{p\in \polyk{2}(\C)\ \vert\ p(0)=c\}$. In particular, the above comment tells us that polynomials in $\poly_{n}^{2,0\to c}(\C)$ is determined uniquely by their critical points.
%\end{definition}
\begin{proof}[Proof of Thm. \ref{thm:LLStratCovBn}]
	The proof is essentially the same in \cite{DM24} but since we do not have to work with equivalence classes, our proof is simpler. We fix a the partitions $\lambda,\mu\vdash n$, with $\lambda\rightarrow\mu$. Suppose the map
	\begin{align*}
		\poly_{\lambda\to\mu}(\C)\to \multi_\mu(\C)
	\end{align*}
	is a local homeomorphism, then the restricted $\LL_B$ is a covering map since $\poly_{\lambda\to\mu}(\C)$ is a surjective proper local homeomorphism from a Hausdorff space to a locally compact Hausdorff space \cite{ho1975note}.
	
	What is left to show is that the restricted $\LL_B$ is a local homeomorphism. This essentially follows from the factorisation of the map Jacobian of the map $\Theta$.
	
	The strata $\C_{[\lambda']}$ can be naturally identified with some open subset of $\C^k$ where $k$ is the number of blocks in $[\lambda']$. 
	Notice that the map $\Theta$ is a local homeomorphism if and only if $z_i\neq 0$ and $z_i\neq z_j$ for all $i\neq j$. Correspondingly, the restriction  $\Theta:\C_{[\lambda']}\to\C^{[\lambda']}$ is a local homeomorphism.
	
	Now consider the open subset $\C_{[\mu']}\subset\C^{[\lambda']}$. The restriction of $\Theta$ to its preimage $\Theta^{-1}(\C_{[\mu']})\subset \C_{[\lambda']}$ is a local homeomorphism. We have the following commuting diagram:
	\begin{displaymath}
		\begin{tikzcd}
			\C_{[\lambda']}\arrow[r,"\Theta"]&\C^{[\lambda']}\\
			\Theta^{-1}(\C_{[\mu']})\arrow[r,"\Theta"]\arrow[u,hook]\arrow[from=2-1, to=4-1]&\C_{[\mu']}\arrow[u,hook]\arrow[d]\\
			&\C_{[\mu]}\arrow[d,"/\sym_\ell"]\\
			\poly_{\lambda\to\mu}(\C)\arrow[r,"\LL_B"]&\multi_\mu(\C)
		\end{tikzcd}
	\end{displaymath}
	The map $\C_{[\mu']}$ to $\C_{[\mu]}$ reduces the shape, but since both spaces are homeomorphic to $\C^\ell$, where $\ell$ is the number of blocks in $[\mu]$, the map is a local homeomorphism. The quotient map by $\sym_\ell$ is shown to be a covering map when restricted to a statum in \cite{DM24} Theorem 6.19.  The map $\C_{[\lambda']}\to\poly_{\lambda}(\C)$ is defined in Proposition~\ref{prop:critPointstoPolyBn}, which is shown to be a covering map. Therefore, after a little diagram chasing, $\LL_B$ is a local homeomorphism.
\end{proof}
\section{A determinant identity}\label{appendix:crit}
In this section, we prove the determinant identity used in the previous appendix.

Consider the function
\begin{align}\label{eqn:critpt}
	p(z)=\int_0^zw^{a_0p+p-1}(w^p-z_1^p)^{a_1}\cdots (w^p-z_m^p)^{a_m}\ dw
\end{align}
where $p,a_i\in\{0,1,2,\cdots\}$. This gives a polynomial of degree $p(1+a_0+\cdots+a_m)$, whose critical points are at $ \{0,z_1,\zeta z_1,\cdots,\zeta^{p-1}z_1,z_2,\cdots,\zeta^{p-1}z_m\} $, where $ \zeta$ is the primitive $ p^{th} $ root of unity. In addition, there are $ m+1 $ critical values corresponding to the values of $ p(0) ,\cdots, p(z_m)$, with multiplicity $ a_0p+p-1,pa_1,\cdots,pa_m $ respectively.

We define an evaluation map $ \theta:\C^m\to\C^m $, which sends $ \mathbf{z}=(z_1,\cdots,z_m) $ to $ \theta(\mathbf{z})=(\theta(z_1),\cdots,\theta(z_m)) $. This map sends the set of nonzero critical points to the set of  corresponding critical values. Let $ \mathbf{J}=\mathbf{J}(\theta) $ be the Jacobian matrix of $\theta(\mathbf{z})$, ie its $(i,j)$-entry is given by $(\mathbf{J})_{ij}=\frac{\D}{\D z_i}p(z_j)$. The study of $\mathbf{J}$ was initiated by Beardon, Carne and Ng \cite{beardon2002critical} who showed that it is invertible for the case $a_0=0$ and $p=a_i=1$. The result was then extended by Dougherty and McCammond \cite{dougherty2020critical} who proved the following:
\begin{thm}\label{thm:DMCrit}
	Let $\mathbf{z}=(z_1,\cdots,z_m)\in \C^m$, $\mathbf{a}=(a_1,\cdots,a_m)$ be an $m$-tuple of positive integers such that $a_1+\cdots+a_m=n$, and let $\mathbf{J}_1$ be the Jacobian defined above with $a_0=0$ and $p=1$. The Jacobian determinant factors as follows:
	\begin{align*}
		\det \mathbf{J_1}=\dfrac{1}{\binom{n}{a_1,\cdots,a_m}}\left(\prod_{j\in [m]}(-z_j)^{a_j}\right)\left(\prod_{\substack{i,j\in[m]\\ j\neq k}}(z_k-z_j)^{a_j}\right)
	\end{align*} 
	Thus, $\mathbf{J}_1$ is invertible if and only if $z_1,\cdots,z_m$ are distinct and nonzero. 
\end{thm}
We first extend their result to the following:
\begin{thm}\label{thm:detidentA2}
	Define the following function:
	\begin{align*}
		p(z)=\int_0^zw^{a_0}(w-z_1)^{a_1}\cdots(w-z_m)^{a_m}\ dw
	\end{align*}
	and $\theta((z_1,\cdots,z_m))=(p(z_1),\cdots,p(z_m))$. Then the determinant of the Jacobian, $\mathbf{J}$, factors as follows:
	\begin{align*}
		\det \mathbf{J}=C\left(\prod_{j\in [m]}(z_j)^{a_0+a_j}\right)\left(\prod_{\substack{i,j\in[m]\\ j\neq k}}(z_k-z_j)^{a_j}\right)
	\end{align*} 
	where $C$ is a nonzero constant.
\end{thm}
%\ben{WHY IS C NONZERO?}
\begin{proof}[Sketch of proof]
	We give a very brief sketch of the proof since the details are largely similar to the proof in \cite{dougherty2020critical}. 
	
	\cite{dougherty2020critical} defined the following polynomial
	\begin{align*}
		q(z_0)=\dfrac{(w-z_0)}{a_0!}^{a_0}\dfrac{(w-z_1)}{a_1!}^{a_1}\cdots\dfrac{(w-z_m)}{a_m!}^{a_m}
	\end{align*}
	We borrow their notation $\tilde{q}_{\mathbf{a},j}$ to be the antiderivative obtained by integrating $q(z_0)$ by parts, integrating the factor $(w-z_j)$ at every step until the other factors become a constant. 
	%we see that $p(z)=\tilde{q}_{\mathbf{a},j}(z)-\tilde{q}_{\mathbf{a},j}(0)$. 
	We obtain the same conclusion as the authors:
	\begin{enumerate}
		\item $\tilde{q}_{\mathbf{a},j}(z_j)=0$
		\item $\tilde{q}_{\mathbf{a},j}(0)$ is divisible by $z_j^{a_j+1}$
		\item $\tilde{q}_{\mathbf{a},j}(z_i)$ is divisible by $(z_i-z_j)^{a_i+a_j+1}$
	\end{enumerate}
	In particular, if we let $z_0=0$, $\tilde{q}_{\mathbf{a},0}(z_i)$ is divisible by $z_i^{a_0+a_i+1}$. Let $\tilde{p}_{\mathbf{a},j}$ be the antiderivative obtained from $\tilde{q}_{\mathbf{a},j}$ by setting $z_0=0$. Notice that $p(z)=\tilde{p}_{\mathbf{a},j}(z)-\tilde{p}_{\mathbf{a},j}(0)$.
	
	The entries of the Jacobian is given by 
	\begin{align*}
		\mathbf{J}_{ij}=a_0!\cdots a_m! \dfrac{\D}{\D z_i}[\tilde{p}_{\mathbf{a},k}(z_j)-\tilde{p}_{\mathbf{a},k}(0)] 
	\end{align*}
	where $k$ can be integer from $\{0,1,2,\cdots,m\}$.
	
	By the divisibility condition of $\tilde{p}_{\mathbf{a},k}$ and choosing appropriate values of $k$, using the arguments in \cite{dougherty2020critical}, one can show that the $i^{th}$ column of $J$ is divisible by $z_i^{a_0+a_i}$ while the difference of $i^{th}$ and $j^{th}$ columns is divisible by $(z_i-z_j)^{a_i+a_j}$. If we rewrite $(z_i-z_j)^{a_i+a_j}=(-1)^{a_i}(z_i-z_j)^{a_j}(z_j-z_i)^{a_i}$, what we have effectively shown is that
	\begin{align*}
		\det\mathbf{J}=C\left(\prod_{j\in [m]}z_j^{a_0+a_j}\right)\left(\prod_{\substack{i,j\in[m]\\ j\neq k}}(z_k-z_j)^{a_j}\right)
	\end{align*}
	The above is an equality because the determinant is a homogenous polynomial of at most $m(a_0+\cdots+a_m)$ degree. Hence, the determinant can only differ from the product by a constant $C$. Finally, setting $a_0=0$, we obtain the result in \cite{dougherty2020critical}. Thus, $C$ must be a nonzero constant.
\end{proof}
\begin{rmk}
	From computer calculations, the constant $C$ should be the multinomial coefficient. 
	\begin{align*}
		C=(-1)^{a_1+\cdots+a_m}\dfrac{1}{\binom{n}{a_0,a_1,\cdots,a_m}}
	\end{align*}
We think the arguments in \cite{dougherty2020critical} can be adapted to show this fact, but since it is not used in this article, we omitted it.
\end{rmk}

\begin{thm}
	Define the following function:
	\begin{align*}
		p(z)=\int_0^zw^{a_0p+p-1}(w^p-z_1^p)^{a_1}\cdots(w^p-z_m^p)^{a_m}\ dw
	\end{align*}
	and $\theta((z_1,\cdots,z_m))=(p(z_1),\cdots,p(z_m))$. Then, the determinant of the Jacobian, $\mathbf{J}$, factors as follows:
	\begin{align*}
		\det \mathbf{J}=C\left(\prod_{j\in [m]}(z_j)^{(a_0+a_j)p+p-1}\right)\left(\prod_{\substack{i,j\in[m]\\ j\neq k}}(z_k-z_j)^{pa_j}\right)
	\end{align*} 
\end{thm}
\begin{proof}
	Let us first rewrite $p(z)$ using the substitution $u=w^p$:
	\begin{align*}
		p(z)&=	\dfrac{1}{p}	\int_0^{z^p}u^{a_0}(u-z_1^p)^{a_1}\cdots (u-z_n^p)^{a_n}\ du\\
	\end{align*}
	This allows us to write $p(z)$ as a composition $q(z^p)$, where $q(z)$ is defined as the function below:
	\begin{align*}
		q(z)&=\dfrac{1}{p}	\int_0^{z}u^{a_0}(u-z_1^p)^{a_1}\cdots (u-z_n^p)^{a_n}\ du
	\end{align*} 
	
	We can now apply chain rule and Theorem~\ref{thm:detidentA2}. Let us denote $\theta_q(z_1,\cdots,z_n)=(q(z_1),\cdots,q(z_n))$ and $\mathbf{z}^p=(z_1^p,\cdots,z_n^p)$, then $\theta$ can be written as $\theta=\theta_q\circ \mathbf{z}^p$
	Thus,
	\begin{align*}
		\mathbf{J}(\theta)&=\left(\mathbf{J}(\theta_q)\circ \mathbf{z}^p\right)\mathbf{J}(\mathbf{z}^p)\\
		&=\left(\mathbf{J}(\theta_q)\circ \mathbf{z}^p\right)\begin{pmatrix}
			pz_1^{p-1}&0&\cdots&0\\
			0&pz_2^{p-1}&0\\
			\vdots&\\
			0&&&pz_n^{p-1}
		\end{pmatrix}
	\end{align*}
	
	Taking the determinant of the Jacobian matrix gives us our desired result:
	\begin{align*}
		\det(\mathbf{J}(\theta))&=\det=\left(\mathbf{J}(\theta_q)\circ \mathbf{z}^p\right)\det\begin{pmatrix}
			pz_1^{p-1}&0&\cdots&0\\
			0&pz_2^{p-1}&0\\
			\vdots&\\
			0&&&pz_n^{p-1}
		\end{pmatrix}\\
		&=C\left(\dfrac{1}{p}\right)^n\left(\prod_{j\in [m]}(z_j^p)^{a_0+a_j}\right)\left(\prod_{\substack{i,j\in[m]\\ j\neq k}}(z_k^p-z_j^p)^{a_j}\right)(p^nz_1^{p-1}\cdots z_n^{p-1})\\
		&=C\left(\prod_{j\in[m]}z_j^{p(a_0+a_j)+p-1}\right)\left(\prod_{\substack{j,k\in[m]\\j\neq k}}(z_k^p-z_j^p)^{a_j}\right)
	\end{align*}
\end{proof}

\end{document}